\documentclass[preprint,12pt]{elsarticle}

\usepackage{amssymb}

\usepackage{amsmath}
  \usepackage{paralist}
  \usepackage{graphics} %% add this and next lines if pictures should be in esp format
  \usepackage{epsfig} %For pictures: screened artwork should be set up with an 85 or 100 line screen
 \usepackage[colorlinks=true]{hyperref}

\allowdisplaybreaks[4]

\usepackage{graphicx,amssymb,mathrsfs,latexsym,amsthm}
\usepackage{verbatim}

\newtheorem{thm}{Theorem}
\newtheorem{lem}[thm]{Lemma}
\newtheorem{cor}[thm]{Corollary}
\newtheorem{prop}[thm]{Proposition}

\newproof{pf}{Proof}

\begin{document}

\begin{frontmatter}

%% Title, authors and addresses

%% use the tnoteref command within \title for footnotes;
%% use the tnotetext command for theassociated footnote;
%% use the fnref command within \author or \address for footnotes;
%% use the fntext command for theassociated footnote;
%% use the corref command within \author for corresponding author footnotes;
%% use the cortext command for theassociated footnote;
%% use the ead command for the email address,
%% and the form \ead[url] for the home page:
%% \title{Title\tnoteref{label1}}
%% \tnotetext[label1]{}
%% \author{Name\corref{cor1}\fnref{label2}}
%% \ead{email address}
%% \ead[url]{home page}
%% \fntext[label2]{}
%% \cortext[cor1]{}
%% \address{Address\fnref{label3}}
%% \fntext[label3]{}

\title{Relative local variational principles for subadditive potentials}

%% use optional labels to link authors explicitly to addresses:
%% \author[label1,label2]{}
%% \address[label1]{}
%% \address[label2]{}

\author[Shanghai]{Xianfeng Ma}
\ead{xianfengma@gmail.com}
\author[Nanjing]{Ercai Chen}
\ead{ecchen@njnu.edu.cn}
\address[Shanghai]{Department of Mathematics, East China University of Science and
Technology\\ Shanghai 200237, China}
\address[Nanjing]{School of Mathematics and Computer Science, Nanjing Normal
University\\Nanjing 210097, China\\
and\\
Center of Nonlinear Science, Nanjing University\\Nanjing 210093,
China}

\begin{abstract}
We prove two relative local variational principles of topological
pressure functions $P(T,\mathcal{F},\mathcal{U},y)$ and
$P(T,\mathcal{F},\mathcal{U}|Y)$ for a given factor map $\pi$, an
open cover $\mathcal{U} $ and a subadditive sequence of real-valued
continuous functions $\mathcal{F}$. By proving the upper
semi-continuity and affinity of the entropy maps
$h_{\{\cdot\}}(T,\mathcal{U}\mid Y)$ and
$h^+_{\{\cdot\}}(T,\mathcal{U}\mid Y)$ on the space of all invariant
Borel probability measures, we show that the relative local pressure
$P(T,\mathcal{\{\cdot\}},\mathcal{U}|Y)$ for subadditive potentials
determines the local measure-theoretic conditional entropies.
\end{abstract}

\begin{keyword}
Pressure\sep variational principle \sep upper semi-continuity \sep
subadditive potentials

\MSC 37A35 \sep 37B40

%% keywords here, in the form: keyword \sep keyword

%% PACS codes here, in the form: \PACS code \sep code

%% MSC codes here, in the form: \MSC code \sep code
%% or \MSC[2008] code \sep code (2000 is the default)

\end{keyword}

\end{frontmatter}

%% \linenumbers

%% main text
\section{Introduction}\label{sec1}
Let $(X,T)$ be a {\it topological dynamical system} (TDS for short)
in the sense that $X$ is a compact metric space and $T:X\rightarrow
X$ is a surjective and continuous map, $\pi$ is a factor map between
TDS $(X,T)$  and $(Y,S)$. The notion of topological pressure was
introduced by Ruelle \cite{Ruelle} for an expansive dynamical system
and later by Walters \cite{Walters1975} for general case. It is
well-known that there exists a basic relationship between the
topological pressure and the relative measure-theoretic entropy.
Ledrappier and Walters \cite{LedWal} formulated the following
classical relative variational principle of pressure for each
$S$-invariant measure $\nu$ on $Y$:
\begin{equation*}
\sup_{\mu\in\mathcal{M}(X,T)}\{h_{\mu}(T,X\mid Y)+\int_Xf(x)d\mu(x):
 \pi\mu=\nu\}=\int_YP(T,f,y)d\nu(y),
\end{equation*}
where $\mathcal{M}(X,T)$ is the family of all $T$-invariant measures
on $X$, $f$ is a real-valued function, $P(T,f,y)$ is the topological
pressure on the compact subset $\pi^{-1}y$, and, for each
$\mu\in\mathcal{M}(X,T)$, $h_{\mu}(T,X\mid Y)$ is the relative
measure-theoretic entropy of $\mu$. For the trivial system $(Y,S)$,
this is the standard variational principle  presented by Walters
\cite{Walters1975}:
\begin{equation*}
\sup_{\mu\in\mathcal{M}(X,T)}\{h_{\mu}(T)+\int_Xf(x)d\mu(x)\}=P(T,f).
\end{equation*}

The topological pressure for nonadditive sequence of potentials has
proved valuable tool in the study of multifractal formalism of
dimension theory, especially for nonconformal dynamical systems
\cite{Cao2008,Barreira2006,Falconer}. Falconer \cite{Falconer} first
introduced the topological pressure for subadditive sequence of
potentials on mixing repellers. He proved the variational principle
for the topological pressure under some Lipschitz conditions and
bounded distortion assumptions on the subadditive potentials. Cao
{\it et al.} \cite{Cao} extended this notion to general compact
dynamical systems, and obtained a subadditive version of variational
principle without any additional assumption. More precisely, let
$\mathcal {F}=\{f_n:n\in\mathbb{N}\}$ be a subadditive sequence of
functions on the TDS, and $\mu(\mathcal{F})=\lim_{n\rightarrow
\infty}\frac{1}{n}\int f_n d\mu$, then
\begin{equation*}
P(T,\mathcal{F})=\sup\{h_{\mu}(T)+\mu(\mathcal{F}):\mu \in
\mathcal{M}(X,T), \mu(\mathcal{F})\neq \infty\}.
\end{equation*}

Since Blanchard \cite{Blanchard1993,Blanchard1995} introduced the
notion of entropy pairs, much attention has been paid to the study
of the local version of the variational principle. Huang {\it et
al.} \cite{huang2007} introduced the notion of local pressure
$P(T,f,\mathcal{U})$, proved the local variational principle of
pressure:
\begin{equation*}
P(T,f,\mathcal{U})=\sup\{h_{\mu}(T,\mathcal{U})+\int_Xf(x)d\mu(x):\mu\in\mathcal{M}(X,T)\},
\end{equation*}
where $h_{\mu}(T,\mathcal{U})$ is the measure-theoretic entropy
relative to $\mathcal{U}$, established the upper semi-continuity and
affinity of the entropy map $h_{\{\cdot\}}(T,\mathcal{U})$, and
showed that the local pressures determine local  measure-theoretic
entropies, i.e., for each $\mu\in \mathcal{M}(X,T)$,
\begin{enumerate}[(a)]
\item
\begin{equation*}
h_{\mu}(T,\mathcal{U})=\inf_{f\in
C(X,\mathbb{R})}\{P(T,f,\mathcal{U})-\int_Xfd\mu\};
\end{equation*}
\item
and if, in addition, $(X,T)$ is invertible, then
\begin{equation*}
h^+_{\mu}(T,\mathcal{U})\leq\inf_{f\in
C(X,\mathbb{R})}\{P(T,f,\mathcal{U})-\int_Xfd\mu\}
\end{equation*}
\end{enumerate}
Zhang \cite{zhang2009} introduced two notions of measure-theoretic
pressure $P_{\mu}^-(T,\mathcal{U},\mathcal {F})$ and
$P_{\mu}^+(T,\mathcal{U},\mathcal {F})$ for a sub-additive sequence
$\mathcal {F}$ of a  real-valued continuous functions on $X$, proved
a local variational principle between topological and
measure-theoretic pressure:
\begin{equation*}
P(T,\mathcal{F},\mathcal{U})=\max_{\mu\in
\mathcal{M}(X,T)}P^-_{\mu}(T,\mathcal{F},\mathcal{U})=\max_{\mu\in
\mathcal{M}(X,T)}\{h_{\mu}(T,\mathcal{U})+\mu(\mathcal{F})\},
\end{equation*}
 and showed the upper semi-continuity of the entropy map
$h^+_{\{\cdot\}}(T,\mathcal{U}) $.

Huang {\it et  al.} \cite{Huang2006} introduced the topological
conditional entropy $h(T,\mathcal{U}\mid Y)$, two notions of
measure-theoretic conditional entropy for covers, i.e.,
$h_{\mu}(T,\mathcal{U}\mid Y)$ and $h_{\mu}^+(T,\mathcal{U}\mid Y)$,
and showed that for a factor map $\pi$ and a given open cover
$\mathcal{U}$, the corresponding variational principles for
conditional entropies hold:
\begin{equation*}
h(T,\mathcal{U}\mid
Y)=\max_{\mu\in\mathcal{M}(X,T)}\{h_{\mu}(T,\mathcal{U}\mid
Y)\},\quad
 h(T,\mathcal{U}\mid
Y)=\max_{\mu\in\mathcal{M}(X,T)}\{h_{\mu}^+(T,\mathcal{U}\mid Y)\}.
\end{equation*}
Zhang  \cite{zhangthesis} introduced the relative local topological
entropy $h(T,\mathcal{U},y)$ and obtained the following relative
local variational principle of the conditional entropy:
\begin{equation*}
\max\{h_{\mu}(T,\mathcal{U}\mid Y):\mu\in\mathcal{M}(X,T) \,\,{\rm
and \,\,}\pi\mu=\nu\}=\int_Yh(T,\mathcal{U},y)d\nu(y).
\end{equation*}

Ma {\it et al.} \cite{Xianfeng2009} and Yan {\it et al.} \cite{Yan}
independently introduced the relative local topological pressure
$P(T,f, \mathcal{U},y)$ for each $y\in Y$. Using the method of
proving the relative variational principle for topological pressure
in \cite{LedWal} and the technique of establishing the conditional
variational principle for the fiber entropy in \cite{Dow},
respectively, they proved the relative local variational principle
for each $\nu\in \mathcal{M}(Y,S)$:
\begin{equation*}
\max_{\mu\in\mathcal{M}(X,T)}\{h_{\mu}(T,\mathcal{U}\mid Y)+ \int_X
f(x) d\mu(x) :\pi\mu=\nu \}=\int_Y P(T,f,\mathcal{U},y)d\nu(y).
\end{equation*}
Yan {\it et al.} \cite{Yan} also proved that the pressure function
$P(T,f, \mathcal{U},y)$ determine the local measure-theoretic
conditional entropy:
\begin{equation*}
h_{\mu}(T,\mathcal{U}|Y)=\inf\{\int_Y
P(T,f,\mathcal{U},y)d\nu(y)-\int_Xfd\mu:f\in C(X,\mathbb{R})\},
\end{equation*}
and obtained the relative local variational principle for the
pressure $P(T,f,\mathcal{U}|Y)$:
\begin{equation*}
P(T,f,\mathcal{U}|Y)=\max_{\mu\in
\mathcal{M}(X,T)}\{h_{\mu}(T,\mathcal{U}|Y)+\int_Xfd\mu\}.
\end{equation*}

The purpose of this paper is to generalize all the  results above to
the case of the relative local topology pressure functions. In fact,
we introduced the two relative local pressure functions $
P(T,\mathcal{F}, \mathcal{U},y)$ and
$P(T,\mathcal{F},\mathcal{U}|Y)$   for subadditive sequence of
potentials, and derive two corresponding relative local variational
principles of pressure. Moreover, we establish the upper
semi-continuity and affinity of the measure-theoretic conditional
entropy maps $h_{\{\cdot\}}(T,\mathcal{U}\mid Y)$ and
$h^+_{\{\cdot\}}(T,\mathcal{U}\mid Y)$, and prove that the relative
local topological pressure $P(T,\mathcal{F},\mathcal{U}|Y)$
determines the measure-theoretic conditional entropies
$h_{\{\cdot\}}(T,\mathcal{U}\mid Y)$ and
$h^+_{\{\cdot\}}(T,\mathcal{U}\mid Y)$. The methods we used is in
the framework of the elegant proof of Huang {\it et al.}
\cite{Huang2006,huang2007} and Ledrappier {\it et al.}
\cite{LedWal}. Our main results state as follows.

\begin{thm}\label{thm4}
Let $\pi:(X,T)\rightarrow (Y,S)$ be a factor map between two TDS and
$\mathcal{U}\in \mathcal{C}_X$. Then the local measure-theoretic
conditional entropy map $h^+_{\{\cdot\}}(T,\mathcal{U}|Y)$ and
$h_{\{\cdot\}}(T,\mathcal{U}|Y)$  are {\it upper semi-continuous and
affine} on $\mathcal{M}(X,T)$.
\end{thm}

\begin{thm}\label{TH}
Let $(X,T)$ be a TDS, $\mathcal{F}\in \mathcal{S}_X$ and
$\mathcal{U}\in \mathcal{C}_X^o$. Let $\pi:(X,T)\rightarrow (Y,S)$
be a factor map between TDS, $\nu\in \mathcal{M}(Y,S)$. Then
\begin{equation*}
\sup_{\mu\in\mathcal{M}(X,T)}\{h_{\mu}(T,\mathcal{U}|Y)
+\mu(\mathcal{F}):\pi\mu=\nu\}=\int_YP(T,\mathcal{F},\mathcal{U},y)d\nu(y).
\end{equation*}
\end{thm}

\begin{thm}\label{coro1}
Let $(X,T)$ be a TDS, $\mathcal{F}\in \mathcal{S}_X$
 and $\mathcal{U}\in \mathcal{C}_X^o$. Let $\pi:(X,T)\rightarrow
(Y,S)$ be a factor map between TDS. Then
\begin{equation*}
\sup\{h_{\mu}(T,\mathcal{U}\mid Y)+\mu(\mathcal{F}): \mu\in
\mathcal{M}(X,T)\}=P(T,\mathcal{F},\mathcal{U}|Y).
\end{equation*}
\end{thm}

\begin{thm}\label{thm2}
Let $(X,T)$ be a TDS, $\mathcal{F}\in \mathcal{S}_X$ and
$\pi:(X,T)\rightarrow (Y,S)$ be a factor map between TDS. Then for
given $\mathcal{U}\in \mathcal{C}_X^o$ and $\mu\in\mathcal{M}(X,T)$,
\begin{equation*}
h_{\mu}(T,\mathcal{U}|Y)=\inf\{P(T,\mathcal{F},\mathcal{U}|Y)-\mu(\mathcal{F}):\mathcal{F}\in
\mathcal{S}_X\}.
\end{equation*}
\end{thm}

\begin{thm}\label{thm3}
Let $(X,T), (Y,S)$ be invertible TDSs, $\mathcal{F}\in
\mathcal{S}_X$, $\pi:(X,T)\rightarrow (Y,S)$ be a factor map between
TDS. Then for given $\mathcal{U}\in\mathcal{C}_X^o$ and $\mu \in
\mathcal{M}(X,T)$,
\begin{equation*}
h^+_{\mu}(T,\mathcal{U}|Y)\leq
\inf\{P(T,\mathcal{F},\mathcal{U}|Y)-\mu(\mathcal{F}):\mathcal{F}\in
\mathcal{S}_X\}.
\end{equation*}
\end{thm}

By Theorem \ref{thm2} and Theorem \ref{thm3}, we immediately obtain
the following result.

\begin{cor}[\cite{Huang2006}]
Let $(X,T), (Y,S)$ be invertible TDSs,  $\pi:(X,T)\rightarrow (Y,S)$
be a factor map between TDS. Then for given
$\mathcal{U}\in\mathcal{C}_X^o$ and $\mu \in \mathcal{M}(X,T)$,
\begin{equation*}
h^+_{\mu}(T,\mathcal{U}|Y)=h_{\mu}(T,\mathcal{U}|Y).
\end{equation*}
\end{cor}

This paper is organized as follows. In Section \ref{sec2}, we
introduce the relative local pressure functions $ P(T,\mathcal{F},
\mathcal{U},y)$ and $P(T,\mathcal{F},\mathcal{U}|Y)$ for subadditive
sequence of potentials and give some necessary lemmas. In Section
\ref{sec3}, we recall some basic properties of the local
measure-theoretic conditional entropies and prove the upper
semi-continuity and affinity of the entropy maps
$h^+_{\{\cdot\}}(T,\mathcal{U}|Y)$ and
$h_{\{\cdot\}}(T,\mathcal{U}|Y)$. In Section \ref{sec4}, we state
and prove the two relative local variational principles for the
topological pressure functions $ P(T,\mathcal{F}, \mathcal{U},y)$
and $P(T,\mathcal{F},\mathcal{U}|Y)$, respectively. In section
\ref{sec5}, using the results we obtained in the former sections, we
prove that the pressure function $P(T,\mathcal{F},\mathcal{U}|Y)$
determines the local measure-theoretic conditional entropies.

\section{Relative local pressure functions for subadditive
potentials}\label{sec2}

Let $(X,T)$ be a TDS and $\mathcal{B}(X)$ be the collection of all
Borel subsets of $X$. Denote by $ \mathcal{M}(X) $ the set of all
Borel, probability measures on $X$, $\mathcal{M}(X,T)$ the set of
$T$-invariant measures, and $\mathcal{M}^e(X,T)$ the set of ergodic
measures. Then
$\mathcal{M}^e(X,T)\subset\mathcal{M}(X,T)\subset\mathcal{M}(X)$,
and $ \mathcal{M}(X), \mathcal{M}(X,T)$ are convex, compact metric
spaces endowed with the weak*-topology. Recall that a {\it cover} of
$X$ is a finite family of Borel subsets of $X$ whose union is $X$,
and, a {\it partition} of $X$ is a cover of $X$ whose elements are
pairwise disjoint. We denote the set of covers, partitions, and open
covers, of $X$, respectively, by $\mathcal{C}_X$, $\mathcal{P}_X$,
$\mathcal{C}_X^o$, respectively. For given two covers $\mathcal{U}$,
$\mathcal{V}\in \mathcal{C}_X$, $\mathcal{U}$ is said to be {\it
finer} than $\mathcal{V}$ (denote by
$\mathcal{U}\succeq\mathcal{V}$) if each element of $\mathcal{U}$ is
contained in some element of $\mathcal{V}$. Let
$\mathcal{U}\vee\mathcal{V}=\{U\cap V: U\in \mathcal{U}, V\in
\mathcal{V}\}$. Given integers $M, N$ with $0\leq M\leq N$ and
$\mathcal{U}\in \mathcal{C}_X$ or $\mathcal{P}_X$, we denote
$\mathcal{U}_M^N$=$\bigvee_{n=M}^NT^{-n}\mathcal{U}$.

Let $(X,T)$  and $(Y,S)$ be two TDS. A continuous map
$\pi:X\rightarrow Y$ is called a {\it factor map} between $(X,T)$
and $(Y,S)$ if it is onto and $\pi T=S\pi$. In this case, we say
that $(X,T)$ is an {\it extension} of $(Y,S)$ or $(Y,S)$ is a factor
of  $(X,T)$.

Let $\pi:(X,T)\rightarrow (Y,S)$ be a factor map between  TDS. Given
$\mathcal{U}\in \mathcal{C}_X$ and $K\subset X$, put
$N(\mathcal{U}\mid K)=\min \{{\rm the\,\, cardinality\,\,
of}\,\,\mathcal{W}: \mathcal{W}\subset\mathcal{U}, \bigcup_{W\in
\mathcal{W}}W\supset K\}$. When $K=X$, we write $N(\mathcal{U}\mid
K)$ simply by $N(\mathcal{U})$. For $y\in Y$, we write
$N(\mathcal{U}\mid y)=N(\mathcal{U}, \pi^{-1} y)$ and
$H(\mathcal{U}\mid y)=\log N(\mathcal{U}\mid y)$. Clearly, if there
is another cover $\mathcal{V}\succeq \mathcal{U}$ then
$H(\mathcal{V}\mid y)\geq H(\mathcal{U}\mid y)$. In fact, for two
covers $\mathcal{U}, \mathcal{V}$ we have $H(\mathcal{U}\vee
\mathcal{V}\mid y)\leq H(\mathcal{U}\mid y)+H(\mathcal{V}\mid y)$.
Let $N(\mathcal{U}\mid Y)=\sup_{y\in Y}N(\mathcal{U}\mid y)$ and
$H(\mathcal{U}\mid Y)=\log N(\mathcal{U}\mid Y)$. Since
$a_n=H(\mathcal{U}_0^{n-1}\mid Y)$ is a non-negative subadditive
sequence, i.e. $a_{n+m}\leq a_n+a_m$, for all $n, m\in \mathbb{N}$,
then the quality
\begin{equation*}
h(T,\mathcal{U}\mid Y)=\lim_{n\rightarrow
\infty}\frac{1}{n}H(\mathcal{U}_0^{n-1}\mid Y)=\inf_{n\geq
1}\frac{1}{n}H(\mathcal{U}_0^{n-1}\mid Y).
\end{equation*}
is well defined, and  called the {\it conditional entropy of }
$\mathcal{U}$ with respect to $(Y,S)$. The {\it topological
conditional entropy of} $(X,T)$ with respect to $(Y,S)$ is defined
 (see \cite{Huang2006}) by
\begin{equation*}
h(T,X\mid Y)=\sup_{\mathcal{U}\in
\mathcal{C}_X^o}h(T,\mathcal{U}\mid Y).
\end{equation*}
If $(Y,S)$ is a trivial system,  this  is  the standard notion of
topological entropy with respect to covers \cite{Walters}.

Let $C(X,\mathbb{R})$ be the Banach space of all continuous,
real-valued functions on $X$ endowed with the supremum norm. Let
$\mathcal{F}=\{f_n:n\in\mathbb{N}\}$ be a sequence of functions in
$C(X,\mathbb{R})$. $\mathcal{F}$ is called {\it subadditive} if for
any $m, n\in \mathbb{N}$ and $x\in X$,
$$
f_{n+m}(x)\leq f_n(x)+f_m (T^n(x)).
$$
Denote by $\mathcal {S}_X$ the set of all subadditive sequences of
functions in $C(X,\mathbb{R})$. In particular, for each $f\in
C(X,\mathbb{R})$, if we set $f_n(x)=\sum_{i=0}^{n-1}f(T^i(x))$, then
$\mathcal{F}=\{f_n:n\in \mathbb{N}\}\in \mathcal{S}_X$. In this
case, for simplicity we write $\mathcal{F}=\{f\}$. For each $c\in
\mathbb{R}$, we let $\{c\}=\{nc:n\in\mathbb{N}\}$. For
$\mathcal{F}=\{f_n:n\in\mathbb{N}\}$,
$\mathcal{G}=\{g_n:n\in\mathbb{N}\}$, and $a,b\in \mathbb{R}$, we
define $a\mathcal{F}+b\mathcal{G}=\{af_n+bg_n:n\in\mathbb{N}\}$ and
$\mathcal{F}=\sup_{n\in \mathbb{N}}\frac{\|f_n\|}{n}$, where
$\|f\|=\sup_{x\in X}f(x)$. Clearly $a\mathcal{F}+b\mathcal{G}\in
\mathcal{S}_X$, and moreover, $(\mathcal{S}_X,\|\cdot\|)$ forms a
Banach space.

If $\nu\in \mathcal{M}(X)$, then for each $n,m \in\mathbb{N}$, $\int
f_{n+m} d\nu \leq \int f_n d\nu +\int f_m d(T^n\nu)$. Thus if
$\mu\in \mathcal{M}(X,T)$, then the sequence $\{\int f_n d\mu : n\in
\mathbb{N}\}$ is subadditive,  so we can set
\begin{equation*}
\mu(\mathcal{F})=\lim_{n\rightarrow \infty }\frac{1}{n}\int
f_n\,d\mu=\inf_{n\in \mathbb{N}}\frac{1}{n}\int f_n\,d\mu \leq
\inf_{n\in \mathbb{N}}\frac{\| f_n\|}{n}.
\end{equation*}
 For each $k\in \mathbb{N}$, let
$\mathcal{F}_k=\{f_{nk}:n\in\mathbb{N}\}$. Then $\mathcal{F}_k$ is a
subsequence of $\mathcal{F}$, and it is easy to see that
$\mathcal{F}_k\in \mathcal{S}_X$ and
$\mu(\mathcal{F}_k)=k\mu(\mathcal{F})$.

For $\mathcal{F}\in \mathcal{S}_X$, $\mathcal{U}\in \mathcal{C}_X^o$
and $y\in Y$, we define
\begin{equation*}
P_n(T,\mathcal{F},\mathcal{U},y)=\inf\{\sum_{V\in\mathcal{V}}\sup_{x\in
V \cap\pi^{-1}(y)}\exp f_n(x): \mathcal{V}\in\mathcal{C}_X \,\, {\rm
and}\,\, \mathcal{V}\succeq \mathcal{U}_0^{n-1}\}.
\end{equation*}
For $V \cap\pi^{-1}(y)=\emptyset$, we let $f_n(x)=-\infty$ for each
$n$. Then the above definition is well defined. Note that for
$\mathcal{F}=\{f\}$, the definition is coincide with that in
\cite{Xianfeng2009}, and for $\mathcal{F}=\{0\}$, it is easy to see
that $P_n(T,\{0\},\mathcal{U},y)=N(\mathcal{U}_0^{n-1},y)$.

For $\mathcal{V}\in \mathcal{C}_X$, we let $\alpha$ be the Borel
partition generated by $\mathcal{V}$ and denote
\begin{equation}
\mathcal{P}^*(\mathcal{V})=\{\beta\in \mathcal{P}_X: \beta\succeq
\mathcal{V} \,\, {\rm and\,\, each \,\, atom \,\, of }\,\,\beta \,\,
{\rm is \,\, the\,\, union\,\, of \,\, some\,\, atoms\,\, of \,\,
}\alpha \}.
\end{equation}

\begin{lem}[\cite{Xianfeng2009}, Lemma 2.1] \label{lem6}
Let $M$ be a compact subset of $X$, $f\in C(X,\mathbb{R})$ and
$\mathcal{V}\in \mathcal{C}_X$. Then
$$
\inf_{\beta\in \mathcal{C}_X, \beta\succeq \mathcal{V}}\sum_{B\in
\beta}\sup_{x\in B\cap M}f(x)=\min\{\sum_{B\in \beta}\sup_{x\in
B\cap M}f(x): \beta\in \mathcal{P}^*(\mathcal{V}) \}.
$$
\end{lem}

If we take $\mathcal{V}=\mathcal{U}_0^{n-1}$, $M=\pi^{-1}(y)$ and
replace $f(x)$ by $\exp f_n(x)$ in Lemma \ref{lem6}, then we have
\begin{equation}\label{eq17}
P_n(T,\mathcal{F},\mathcal{U},y)=\min\{\sum_{B\in\beta}\sup_{x\in
B\cap\pi^{-1}(y)}\exp f_n(x):\beta\in
\mathcal{P}^*(\mathcal{U}_0^{n-1}) \}.
\end{equation}
In particular, if $\mathcal{U}$ is a partition, then
$P_n(T,\mathcal{F},\mathcal{U},y)=\sum_{U\in\mathcal{U}_0^{n-1}}\sup_{x\in
U\cap \pi^{-1}(y)}\exp f_n(x)$.

\begin{lem}\label{lem2}
Let $(X,T)$ be a TDS, $\mathcal{F}\in \mathcal{S}_X$ and
$\mathcal{U}=\{U_1,\cdots,U_d\}\in \mathcal{C}_X^o$. Let
$\pi:(X,T)\rightarrow (Y,S)$ be a factor map between TDS. Then the
mappings $y\rightarrow P_n(T,\mathcal{F},\mathcal{U},y)$ of $Y$ to
$\mathbb{R}$ are universally measurable for any $n\geq 1$ and there
exists a constant $M$ such that $ \frac{1}{n}\log
P_n(T,\mathcal{F},\mathcal{U},y)\leq M $ for all $n\geq 1$ and $y\in
Y$.
\end{lem}

\begin{pf}
The proof of the measurability can be seen in \cite{Xianfeng2009}.
For the other part, since
 $$P_n(T,\mathcal{F},\mathcal{U},y)\leq e^{\|f_n\|}\cdot\min_{\beta\in\mathcal {P}^*(\mathcal{U}_0^{n-1})}{\rm
card}(\beta)\leq e^{\|f_n\|}\cdot d^n.$$ Then $ \frac{1}{n}\log
P_n(T,\mathcal{F},\mathcal{U},y)\leq \frac{\|f_n\|}{n}+\log d. $ Let
$M=\|f_1\|+\log d$, and we  get the result.
\end{pf}

For each $y\in Y$, $\mathcal{U}\in \mathcal{C}_X^o$ and
$\mathcal{F}\in \mathcal{S}_X$, we define the universally measurable
map $P(T,\mathcal{F},\mathcal{U},y)$ from $Y$ to $\mathbb{R}$ as
\begin{equation*}
P(T,\mathcal{F},\mathcal{U},y)=\limsup_{n\rightarrow
\infty}\frac{1}{n}\log P_n (T,\mathcal{F},\mathcal{U},y).
\end{equation*}

For each $\nu \in \mathcal {M}(Y,S)$, the following lemma shows that
the  limit superior in the above definition can be obtained by the
limit for $\nu$-{\it a.e.} $y\in Y$.

\begin{lem} \label{lem5}
 Let $\nu \in \mathcal {M}(Y,S)$. For
$\mathcal{F}\in \mathcal{S}_X$, $\mathcal{U}\in \mathcal{C}_X^o$,
and $\nu$-a.e. $y\in Y$,
$$
P(T,\mathcal{F},\mathcal{U},y)= \lim_{n\rightarrow
\infty}\frac{1}{n}\log P_n (T,\mathcal{F},\mathcal{U},y)
$$
exists.
\end{lem}

\begin{pf}
 For any $n,m\in\mathbb{N}$, $\mathcal{V}_1\succeq
\mathcal{U}_0^{n-1}$, $\mathcal{V}_2\succeq \mathcal{U}_0^{m-1}$, we
have $\mathcal{V}_1\vee T^{-n}\mathcal{V}_2\succeq
\mathcal{U}_0^{n+m-1}$. It follows that
\begin{align*}
P_{n+m}(T,\mathcal{F},\mathcal{U},y)&\leq \sum_{V_1\in
\mathcal{V}_1}\sum_{V_1\in \mathcal{V}_2}\sup_{x\in V_1\cap
T^{-n}V_2\cap \pi^{-1}(y)}\exp f_{n+m}(x)\\
&\leq\sum_{V_1\in \mathcal{V}_1}\sum_{V_2\in
\mathcal{V}_2}\sup_{x\in
V_1\cap T^{-n}V_2\cap \pi^{-1}(y)}\exp(f_n(x)+f_m(T^nx))\\
&\leq \sum_{V_1\in \mathcal{V}_1}\sum_{V_2\in \mathcal{V}_2}
\big(\sup_{x\in V_1\cap \pi^{-1}(y)}\exp f_n(x)\cdot\sup_{z\in
V_2\cap
\pi^{-1}(S^ny)}\exp f_m(z)\big)\\
&=\big(\sum_{V_1\in \mathcal{V}_1}\sup_{x\in V_1\cap
\pi^{-1}(y)}\exp f_n(x)  \big)\big(\sum_{V_2\in
\mathcal{V}_2}\sup_{z\in V_2\cap \pi^{-1}(S^ny)}\exp f_m(z)  \big).
\end{align*}
Since $\mathcal{V}_i,i=1,2$ is arbitrary, then
$P_{n+m}(T,\mathcal{F},\mathcal{U},y)\leq
P_n(T,\mathcal{F},\mathcal{U},y)P_m(T,\mathcal{F},\mathcal{U},S^ny)$,
and so $\log P_n(T,\mathcal{F},\mathcal{U},y) $ is subadditive. By
Kingman's subadditive ergodic theorem (See \cite{Walters}), we
complete the proof.
\end{pf}

The following Lemma  follows from Lemma \ref{lem5} directly.

\begin{lem}\label{lem13}
Let $\nu \in \mathcal {M}(Y,S)$. Then
$P(T^k,\mathcal{F}_k,\mathcal{U}_0^{n-1},y)=kP(T,\mathcal{F},\mathcal{U},y)$
for $\mathcal{F}\in \mathcal{S}_X$,
$\mathcal{U}=\{U_1,\cdots,U_d\}\in \mathcal{C}_X^o$, $k\in
\mathbb{N}$ and $\nu$-a.e. $y\in Y$.
\end{lem}

We refer to $P(T,\mathcal{F},\mathcal{U},y)$ as the {\it topological
pressure of $\mathcal{F}$ relative to $\mathcal{U}$  on $\pi^{-1}y$
}.

Let $$P_n(T,\mathcal{F},\mathcal{U},Y)=\sup_{y\in
Y}P_n(T,\mathcal{F},\mathcal{U},y).$$

\begin{lem}\label{lem16}
For each $\mathcal{U}\in \mathcal{C}_X$, $n\in \mathbb{N}$, there
exists $\eta \in \mathcal{P}_X$ with $\eta \succeq
\mathcal{U}_0^{n-1}$ such that for each $y\in Y$,
\begin{equation*}
\sum_{C\in \eta \cap \pi^{-1}(y)}\sup_{x\in C}(\exp f_n(x))\leq
P_n(T,\mathcal{F},\mathcal{U},Y).
\end{equation*}
\end{lem}

\begin{pf}
For each $y \in Y$, by Lemma \ref{lem6}, there exists $\beta_y\in
\mathcal{P}^*(\mathcal{U}_0^{n-1})$ such that
\begin{equation*}
\sum_{C\in\beta_y\cap\pi^{-1}(y)}\sup_{x\in C}(\exp
f_n(x))=P_n(T,\mathcal{F},\mathcal{U},y)\leq
P_n(T,\mathcal{F},\mathcal{U},Y).
\end{equation*}

Since $\mathcal{P}^*(\mathcal{U}_0^{n-1})$ is finite, we can find
$y_1,y_2,\cdots,y_s\in Y$ such that for each $y\in Y$, there exists
$i\in\{1,2,\cdots,s\}$ such that
$\sum_{C\in\beta_{y_i}\cap\pi^{-1}(y)}\sup_{x\in C}(\exp
f_n(x))=P_n(T,\mathcal{F},\mathcal{U},y)$. For each
$i\in\{1,2,\cdots,s\}$, define
$$D_i=\{ y\in Y:  \sum_{C\in\beta_{y_i}\cap\pi^{-1}(y)}\sup_{x\in C}(\exp
f_n(x))=P_n(T,\mathcal{F},\mathcal{U},y)\}.$$ Let
$C_i=D_i\backslash\bigcup_{j=1}^{i-1}D_j, i=1,2,\cdots,s.$ Then
$C_i\cap C_j=\emptyset, i\neq j$, and it is easy to see that
\begin{equation*}
\eta=\{ \beta_{y_i}\cap \pi^{-1}(C_i): i=1,2,\cdots,s \},
\end{equation*}
where $\beta_{y_i}\cap \pi^{-1}(C_i)=\{B\cap\pi^{-1}(C_i): B\in
\beta_{y_i} \}$, is a partition of $X$ finer than
$\mathcal{U}_0^{n-1}$. Moreover, for each $y\in Y$, there exists
$i\in \{1,2,\cdots,s\}$ such that
\begin{equation*}
\sum_{C\in \eta\cap \pi^{-1}(y)}\sup_{x\in C}(\exp
f_n(x))=\sum_{C\in \beta_{y_i}\cap \pi^{-1}(y)}\sup_{x\in C}(\exp
f_n(x))\leq P_n(T,\mathcal{F},\mathcal{U},Y),
\end{equation*}
and we complete the proof.
\end{pf}

From the proof of Lemma \ref{lem5},   it is not hard to see that the
sequence of functions $\log P_n(T,\mathcal{F},\mathcal{U},Y)$ is
subadditive. The {\it topological pressure of $\mathcal{F}$ relative
to $\mathcal{U}$ and $(Y,S)$} is defined by
\begin{equation*}
P(T,\mathcal{F},\mathcal{U}\mid Y )=\lim_{n\rightarrow
\infty}\frac{1}{n}\log P_n(T,\mathcal{F},\mathcal{U},Y)=\inf_{n\in
\mathbb{N}}\frac{1}{n}\log P_n(T,\mathcal{F},\mathcal{U},Y)
\end{equation*}

The {\it topological pressure of $\mathcal{F}$} is defined by
\begin{equation*}
P(T,\mathcal{F}\mid Y)=\sup_{\mathcal{U}\in
\mathcal{C}_X^o}P(T,\mathcal{F},\mathcal{U}\mid Y).
\end{equation*}

For the trivial system $(Y,S)$, it is not hard to see that the
topological pressure defined $P(T,\mathcal{F}\mid Y)$ above is
equivalent to the ones defined in \cite{zhang2009}. Moreover, if
$(Y,S)$ is the trivial system and $\mathcal{F}=\{f\}$,  then
$P(T,\mathcal{F}, \mathcal{U}\mid Y)$ is the definition defined in
\cite{huang2007}. If $\mathcal{F}=\{0\}$, then
$P(T,\{0\},\mathcal{U}\mid Y)=h(T,\mathcal{U}\mid Y)$. If $(Y,S)$ is
the trivial system and $\mathcal{F}=\{0\}$, then
$P(T,\{0\},\mathcal{U}\mid Y)=h(T,\mathcal{U})$, which is the
standard topological entropy with respect to the cover
$\mathcal{U}$. As in \cite{huang2007}, the advantage of the above
definition of $P_n(T,\mathcal{F},\mathcal{U},y)$ is the
monotonicity, i.e., if $\mathcal{U}\succeq \mathcal{V}$, then
$P_n(T,\mathcal{F},\mathcal{U},y)\geq
P_n(T,\mathcal{F},\mathcal{V},y)$.

\section{Measure-theoretic conditional entropies}\label{sec3}

Given a partition $\alpha\in \mathcal{P}(X)$, $\mu\in
\mathcal{M}(X)$ and a sub-$\sigma$-algebra
$\mathcal{A}\subset\mathcal{B}(X)$, define
\begin{equation*}
H_{\mu}(\alpha\mid\mathcal{A})=\sum_{A\in\alpha}\int_X-\mathbb{E}(1_A\mid\mathcal{A})\log
\mathbb{E}(1_A\mid\mathcal{A})d\mu,
\end{equation*}
where $\mathbb{E}(1_A\mid\mathcal{A})$ is the expectation of $1_A$
with respect to $\mathcal{A}$. One standard fact states that
$H_{\mu}(\alpha\mid\mathcal{A})$ increases with respect to $\alpha$
and decreases with respect to $\mathcal{A}$.

When $\mu \in \mathcal{M}(X,T)$ and $\mathcal{A}$ is a $T$-invariant
$\mu$-measurable $\sigma$-algebra of $X$, i.e.
$T^{-1}\mathcal{A}\subset\mathcal{A}$,
$H_{\mu}(\alpha_0^{n-1}\mid\mathcal{A})$ is a non-negative
subadditive sequence for a given $\alpha\in \mathcal{P}_X$. The {\it
measure-theoretic conditional entropy of $\alpha$ with respect to
$\mathcal{A}$ } is defined as
\begin{equation}\label{def1}
h_{\mu}(T,\alpha\mid\mathcal{A})
=\lim_{n\rightarrow\infty}\frac{1}{n}H_{\mu}(\alpha_0^{n-1}\mid\mathcal{A})
=\inf_{n\geq 1}H_{\mu}(\alpha_0^{n-1}\mid\mathcal{A}),
\end{equation}
and the {\it measure-theoretic conditional entropy of $(X,T)$ with
respect to $\mu$ } is defined as
\begin{equation*}
h_{\mu}(T,X\mid\mathcal{A})=\sup_{\alpha\in\mathcal{P}_X}h_{\mu}(T,\alpha\mid\mathcal{A}).
\end{equation*}
Particularly, if $\pi: (X,T)\rightarrow (Y,S)$ is a factor map
between TDS and $\alpha\in \mathcal{P}_X$, the {\it conditional
entropy of $\alpha$ with respect to $(Y,S)$ } is defined as
\begin{equation*}
h_{\mu}(T,\alpha\mid Y)=h_{\mu}(T,\alpha\mid
\pi^{-1}(\mathcal{B}(Y)))
=\lim_{n\rightarrow\infty}\frac{1}{n}H_{\mu}(\alpha_0^{n-1}\mid\pi^{-1}(\mathcal{B}(Y))),
\end{equation*}
and the {\it measure-theoretic conditional entropy of $(X,T)$ with
respect to $(Y,S)$ } is defined as
\begin{equation*}
h_{\mu}(T,X\mid Y)=\sup_{\alpha\in\mathcal{P}_X}h_{\mu}(T,\alpha\mid
Y).
\end{equation*}
For the classical theory of measure-theoretic entropy, we refer the
reader to \cite{parry1981, Walters, Ye2008}.

\begin{lem}[\cite{Huang2006}, Lemma 3.3]\label{lem12}
Let $\pi:(X,T)\rightarrow (Y,S)$ be a factor map between two TDS and
$\alpha\in \mathcal{P}_X$. Then the following hold:
\begin{enumerate}
\item The function $H_{\{\cdot\}}(\alpha\mid Y)$ is concave on
$\mathcal{M}(X)$;
\item The function $h_{\{\cdot\}}(\alpha\mid Y)$ and $h_{\{\cdot\}}(T, X\mid
Y)$ are affine on $\mathcal{M}(X,T)$.
\end{enumerate}
\end{lem}

A real-valued function $f$ defined on a compact metric space $Z$ is
called {\it upper semi-continuous }(for short u.s.c.) if one of the
following equivalent conditions holds:
\begin{enumerate}
\item $\limsup_{z'\rightarrow z}f(z')\leq f(z)$ for each $z\in Z$;
\item for each $f\in C(Z, \mathbb{R})$ the set $\{z\in Z:f(z)\geq r\}$ is
closed.\label{usc2}
\end{enumerate}
By \ref{usc2}, the infimum of any family of u.s.c. functions is
again a u.s.c. one; both the sum and supremum of finitely many
u.s.c. functions are u.s.c. ones.

A subset $A$ of $X$ is called {\it clopen} if it is both closed and
open in $X$. A partition is called {\it clopen} if it consists of
clopen sets.

\begin{lem}[\cite{Huang2006}, Lemma 3.4]\label{lem11}
Let $\pi:(X,T)\rightarrow (Y,S)$ be a factor map between two TDS and
$\alpha\in \mathcal{P}_X$ whose elements are clopen sets of $X$.
Then:
\begin{enumerate}
\item  $H_{\{\cdot\}}(\alpha\mid Y)$ is a {\it u.s.c.} function on
$\mathcal{M}(X)$;
\item  $h_{\{\cdot\}}(\alpha\mid Y)$ is a {\it u.s.c.} function on $\mathcal{M}(X,T)$.
\end{enumerate}
\end{lem}

Inspired by the ideas of Romagnoli \cite{Rom2003} in local entropy
for covers, Huang {\it et al.} \cite{Huang2006} introduced a new
notion of $\mu$-measure-theoretic conditional  entropy for covers,
which extends definition \eqref{def1} to covers. Let
$\pi:(X,T)\rightarrow (Y,S)$ be a factor map  and $\mu \in
\mathcal{M}(X)$. For $\mathcal{U}\in \mathcal{C}_X$ define
\begin{equation}
H_{\mu}(\mathcal{U}\mid Y)=\inf_{\alpha\in \mathcal{P}_X,
\alpha\succeq\mathcal{U}}H_{\mu}(\alpha\mid \pi^{-1}\mathcal{B}(Y)).
\end{equation}
In particular, $H_{\mu}(\alpha\mid Y)=H_{\mu}(\alpha\mid
\pi^{-1}\mathcal{B}(Y))$ for $\alpha\in \mathcal{P}_X$. Many
properties of the conditional function $H_{\mu}(\alpha\mid Y)$ for a
partition $\alpha$ can be extended to $H_{\mu}(\mathcal{U}\mid Y) $
for a cover $\mathcal{U}$; for details see \cite{Huang2006}.

\begin{lem}\label{lem7}
Let $\pi:(X,T)\rightarrow (Y,S)$ be a factor map between TDS and
$\mu\in \mathcal{M}(X)$. If $\mathcal{U}, \mathcal{V}\in
\mathcal{C}_X$, then the following hold:
\begin{enumerate}
\item $0\leq H_{\mu}(\mathcal{U}\mid Y) \leq \log N(\mathcal{U})$;
\item if $\mathcal{U}\succeq \mathcal{V}$, then $H_{\mu}(\mathcal{U}\mid Y)\geq H_{\mu}(\mathcal{V}\mid
Y)$;
\item $H_{\mu}(\mathcal{U}\vee\mathcal{V}\mid Y)\leq H_{\mu}(\mathcal{U}\mid
Y)+H_{\mu}(\mathcal{V}\mid Y)$;
\item $H_{\mu}(T^{-1}\mathcal{U}\mid Y)\leq H_{T\mu}(\mathcal{U}\mid
Y)$.
\end{enumerate}
\end{lem}

\begin{lem}[\cite{zhangthesis}, Lemma 5.2.8]\label{lem15}
Let $\pi:(X,T)\rightarrow (Y,S)$ be a factor map between two TDS,
$\mathcal{U}\in \mathcal{C}_X$, $\mu\in \mathcal{M}(X,T)$. Let
$\mu=\int_Y\mu_yd\nu(y)$ be the disintegration of $\mu$ over $\nu$
where $\nu=\pi\mu$. Then
$$
H_{\mu}(\mathcal{U}\mid Y)=\int_Y H_{\mu_y}(\mathcal{U})d\nu(y),
$$
where $H_{\mu_y}(\mathcal{U})=\inf_{\alpha\in \mathcal {P}_X,\alpha
\succeq\mathcal {U}}H_{\mu_y}(\alpha)$.
\end{lem}

For a given $\mathcal{U}\in \mathcal{C}_X$, $\mu \in
\mathcal{M}(X,T)$, it follows easily from Lemma  \ref{lem7} that
$H_{\mu}(\mathcal{U}^{n-1}_0\mid Y)$ is a subadditive function of
$n\in \mathbb{N}$. Hence the local $\mu$-conditional entropy of
$\mathcal{U}$ with respect to $(Y,S)$ can be defined as
\begin{equation}
h_{\mu}(T,\mathcal{U}\mid Y)=\lim_{n\rightarrow
\infty}\frac{1}{n}H_{\mu}(\mathcal{U}_0^{n-1}\mid Y)=\inf_{n\geq
1}\frac{1}{n}H_{\mu}(\mathcal{U}_0^{n-1}\mid Y).
\end{equation}

This extension of local measure-theoretic conditional entropy from
partitions to covers allows the generalization of the relative local
variational principle of entropy to the relative variational
principle of pressure.

Following the works of Romagnoli \cite{Rom2003}, Huang {\it et al.}
\cite{Huang2006} also introduced another type of local
$\mu$-conditional entropy. Let $\pi:(X,T)\rightarrow (Y,S)$ be a
factor map between TDS. Given $\mu\in \mathcal{M}(X,T)$ and
$\mathcal{U}\in \mathcal{C}_X$ define
\begin{equation}
h_{\mu}^+(T,\mathcal{U}\mid Y)=\inf_{\alpha\in
\mathcal{P}_X,\alpha\succeq \mathcal{U}}h_{\mu}(T,\alpha\mid Y).
\end{equation}
Clearly, $h_{\mu}^+(T,\mathcal{U}\mid Y)\geq
h_{\mu}(T,\mathcal{U}\mid Y)$. Moreover, for a factor map between
TDS, the following lemma holds.

\begin{lem}[\cite{Huang2006}, Lemma 4.1(3)]\label{lem10}
Let $\pi:(X,T)\rightarrow (Y,S)$ be a factor map between two TDS and
$\mu\in \mathcal{M}(X,T)$. Then for each $\mathcal{U}\in
\mathcal{C}_X$,
$$
h_{\mu}(T,\mathcal{U}\mid
Y)=\lim_{n\rightarrow\infty}\frac{1}{n}h_{\mu}^+(T^n,\mathcal{U}_0^{n-1}\mid
Y)=\inf_{n\in
\mathbb{N}}\frac{1}{n}h_{\mu}^+(T^n,\mathcal{U}_0^{n-1}\mid Y).
$$
\end{lem}

For each $\mu \in \mathcal{M}(X,T)$, there exists a unique Borel
probability measure $m$ on $\mathcal{M}^e(X,T)$ such that $\mu
=\int_{\mathcal{M}^e(X,T)}\theta dm(\theta)$, i.e. $\mu $ admits an
ergodic decomposition. The ergodic decomposition of $\mu$ gives rise
to an ergodic decomposition of the $\mu$-entropy relative to the
partition $\alpha\in \mathcal{P}_X$:
\begin{equation*}
h_{\mu}(T,\alpha)=\int_{\mathcal{M}^e(X,T)}h_{\theta}(T,\alpha)dm(\theta).
\end{equation*}
Following the ideas of proving the ergodic decompositions of the
$\mu$-entropies  relative to covers \cite{huangye2006}, Huang {\it
et al.} \cite{Huang2006} gave the ergodic decompositions of the two
kinds of measure conditional entropy of covers.

\begin{lem}[\cite{Huang2006}, Lemma 5.3]\label{lem8}
Let $\pi:(X,T)\rightarrow (Y,S)$ be a factor map between two TDS,
$\mu\in \mathcal{M}(X,T)$ and $\mathcal{U}\in \mathcal{C}_X$. If
$\mu=\int_{\mathcal{M}^e(X,T)}\theta dm(\theta)$ is the ergodic
decomposition of $\mu$, then
\begin{enumerate}
\item $h_{\mu}^+(T,\mathcal{U}\mid Y)=\int_{\mathcal{M}^e(X,T)}h_{\theta}^+(T,\mathcal{U}\mid
Y)dm(\theta)$;
\item $h_{\mu}(T,\mathcal{U}\mid Y)=\int_{\mathcal{M}^e(X,T)}h_{\theta}(T,\mathcal{U}\mid
Y)dm(\theta)$.
\end{enumerate}
\end{lem}

Now we are ready  to prove Theorem \ref{thm4}, i.e., if we let
$\pi:(X,T)\rightarrow (Y,S)$ be a factor map between two TDS and
$\mathcal{U}\in \mathcal{C}_X$, then the local conditional entropy
map $h^+_{\{\cdot\}}(T,\mathcal{U}|Y)$ and
$h_{\{\cdot\}}(T,\mathcal{U}|Y)$  are {\it upper semi-continuous and
affine} on $\mathcal{M}(X,T)$.

\begin{pf}[Proof of Theorem \ref{thm4}]

We first prove the upper semi-continuity. Let
$\mathcal{U}=\{U_1,\cdots, U_M\}$. By Lemma \ref{lem10},
$h_{\mu}(T,\mathcal{U}\mid Y)=\inf_{n\in
\mathbb{N}}\frac{1}{n}h_{\mu}^+(T^n,\mathcal{U}_0^{n-1}\mid Y).$ It
follows that if the local conditional entropy map
$h^+_{\{\cdot\}}(T,\mathcal{U}|Y): \mu \in
\mathcal{M}(X,T)\rightarrow \mathbb{R}$ is upper semi-continuous,
then $h_{\{\cdot\}}(T,\mathcal{U}|Y): \mu \in
\mathcal{M}(X,T)\rightarrow \mathbb{R}$ is also upper
semi-continuous.

We now prove $h^+_{\{\cdot\}}(T,\mathcal{U}|Y): \mu \in
\mathcal{M}(X,T)\rightarrow \mathbb{R}$ is upper semi-continuous.
Since for each $\mu \in \mathcal{M}(X,T)$,
\begin{equation*}
h^+_{\mu}(T,\mathcal{U}|Y)=\inf_{\alpha\in \mathcal{P}_X,
\alpha\succeq \mathcal{U}}\inf_{n\in
\mathbb{N}}\frac{1}{n}H_{\mu}(\alpha_0^{n-1}|Y)=\inf_{n\in
\mathbb{N}}\inf_{\alpha\in \mathcal{P}_X, \alpha\succeq
\mathcal{U}}\frac{1}{n}H_{\mu}(\alpha_0^{n-1}|Y),
\end{equation*}
it is suffice to prove that for each $n\in \mathbb{N}$,   the map
$\phi_n(\mu)=\inf_{\alpha\in \mathcal{P}_X, \alpha\succeq
\mathcal{U}}H_{\mu}(\alpha_0^{n-1}|Y)$ is upper semi-continuous on
$\mathcal{M}(X,T)$. Moreover, By the definition of the upper
semi-continuous function, it is suffice to prove that for each $\mu
\in \mathcal{M}(X,T)$ and $\epsilon>0$,
\begin{equation*}
\limsup_{\mu'\rightarrow \mu,\mu \in
\mathcal{M}(X,T)}\phi_n(\mu')\leq \phi_n(\mu)+\epsilon.
\end{equation*}

Fix $\mu \in \mathcal{M}(X,T)$ and $\epsilon >0$. There exists
$\alpha \in \mathcal{P}_X, \alpha\succeq\mathcal{U}$ such that
\begin{equation*}
H_{\mu}(\alpha_0^{n-1}|Y)\leq \phi_n(\mu)+\epsilon/2.
\end{equation*}
Without loss of the generality, we assume that $\alpha=\{A_1,
\cdots, A_M\}$ with $A_i\subset U_i$ for each $1\leq i \leq M$. Let
$\mu^n=\sum_{i=0}^{n-1}T^i\mu$. By Lemma 4.15 \cite{Walters}, there
exists a  $\delta=\delta(M,n,\epsilon)>0$ such that whenever
$\beta^1=\{B_1^1,B_2^1,\cdots,B_k^1\}$ and
$\beta^2=\{B_1^2,B_2^2,\cdots,B_k^2\}$ are k-measurable partitions
with $\sum_{i=1}^k\mu^n(B_i^1 \Delta B_i^2 )<\delta$, then
\begin{equation*}
H_{\mu}(\bigvee_{i=0}^{n-1}T^{-i}\beta^1|\bigvee_{i=0}^{n-1}T^{-i}
\beta^2)\leq \sum_{i=0}^{n-1}H_{T^i\mu}(\beta^1|\beta^2)\leq
H_{\sum_{i=0}^{n-1}T^i\mu}(\beta^1|\beta^2)<\epsilon/2.
\end{equation*}

Let $\mathcal{U}_{\mu,n}^*=\{\beta\in \mathcal{P}_X:\beta\succeq
\mathcal{U} \,\,\text{and} \,\,
\mu(\bigcup_{C\in\beta_0^{n-1}}\partial C)=0\}$. Then there exists
$\beta=\{B_1,\cdots,B_M\}\in \mathcal{U}_{\mu,n}^*$ such that
$\sum_{i=1}^M\mu^n(A_i\Delta B_i)<\delta$ and
$H_{\mu}(\beta_0^{n-1}|\alpha_0^{n-1})<\epsilon/2$ (See Claim P.164
\cite{Ye2008}). Note that the condition
$\mu(\bigcup_{C\in\beta_0^{n-1}}\partial C)=0$ in the definition of
$\mathcal{U}_{\mu,n}^*$ implies that $\mu(\sum_{i=0}^M \partial
B_i)=0$. Then, by Lemma 3.2 (ii) \cite{LedWal},
\begin{align*}
\limsup_{\mu'\rightarrow \mu,\mu'\in \mathcal{M}(X,T)}\phi_n(\mu')
&\leq \limsup_{\mu'\rightarrow \mu,\mu'\in
\mathcal{M}(X,T)}H_{\mu'}(\beta_0^{n-1}|Y) \\
&\leq H_{\mu}(\beta_0^{n-1}|Y)\\
&\leq
H_{\mu}(\alpha_0^{n-1}|Y)+H_{\mu}(\beta_0^{n-1}|\alpha_0^{n-1}\vee
Y)\\
&\leq
H_{\mu}(\alpha_0^{n-1}|Y)+H_{\mu}(\beta_0^{n-1}|\alpha_0^{n-1})\\
&\leq \phi_n(\mu)+\epsilon.
\end{align*}

We now prove the affinity. Given $\mu_i\in \mathcal{M}(X,T)$,
$i=1,2$, and $0<\lambda <1$. Let
$\mu_i=\int_{\mathcal{M}^e}(X,T)\theta dm_i(\theta)$ be the ergodic
decomposition of $\mu_i$. Let $\mu=\lambda\mu_1+(1-\lambda)\mu_2$
and $m=\lambda m_1+(1-\lambda)m_2$. It is clear that $m$ is a Borel
probability measure on $\mathcal{M}^e(X,T)$ and
$\mu=\int_{\mathcal{M}^e}(X,T)\theta dm(\theta)$. By Lemma
\ref{lem8},
\begin{align*}
h^+_{\mu}&(T,\mathcal{U}|Y)=\int_{\mathcal{M}^e(X,T)}h_{\theta}^+(T,\mathcal{U}\mid
Y)dm(\theta)\\
&=\lambda \int_{\mathcal{M}^e(X,T)}h_{\theta}^+(T,\mathcal{U}\mid
Y)dm_1(\theta)
+(1-\lambda)\int_{\mathcal{M}^e(X,T)}h_{\theta}^+(T,\mathcal{U}\mid
Y)dm_2(\theta)\\
&=\lambda
h^+_{\mu_1}(T,\mathcal{U}|Y)+(1-\lambda)h^+_{\mu_2}(T,\mathcal{U}|Y).
\end{align*}
Then the local conditional entropy map
$h^+_{\{\cdot\}}(T,\mathcal{U}|Y)$ is affine on  $\mathcal{M}(X,T)$.
The proof the affinity of $h_{\{\cdot\}}(T,\mathcal{U}|Y)$ is
similar to the above proof.
\end{pf}

For the trivial system $(Y,S)$, it is clear that the following
result holds, which was proved in \cite{huang2007} and
\cite{zhang2009}.

\begin{cor}
Let (X,T) be a TDS and $\mathcal{U}\in\mathcal{C}_X^o$. Then the
local entropy maps $h^+_{\{\cdot\}}(T,\mathcal{U})$ and
$h_{\{\cdot\}}(T,\mathcal{U})$ are upper semi-continuous and affine
on $\mathcal{M}(X,T)$.
\end{cor}

\section{Relative local  variational principles for subadditive potentials}\label{sec4}

\begin{lem}[\cite{zhangthesis}, Proposition 5.2.9]\label{lem9}
Let $\pi:(X,T)\rightarrow (Y,S)$ and $\varphi : (Z,R)\rightarrow
(X,T)$ be two  factor maps between TDS. If $\tau\in
\mathcal{M}(Z,R)$, $\mu=\varphi\tau\in\mathcal{M}(X,T)$, then for
each $\mathcal{U}\in \mathcal{C}_X$,
\begin{equation*}
h_{\tau}(R,\varphi^{-1}(\mathcal{U})\mid
Y)=h_{\mu}(T,\mathcal{U}\mid Y).
\end{equation*}
\end{lem}

\begin{lem}\label{lem4}
Let $\pi:(X,T)\rightarrow (Y,S)$ and $\varphi : (Z,R)\rightarrow
(X,T)$ be two  factor maps between TDS, $\mathcal{F}\in
\mathcal{S}_X$  and $\mathcal{U}\in \mathcal{C}_X^o$. Then for each
$y\in Y$ and $n\in \mathbb{N}$, $P_n(R,\mathcal{F}\circ
\varphi,\varphi^{-1}\mathcal{U},y)=P_n(T,\mathcal{F},\mathcal{U},y)$,
where $\mathcal{F}\circ \varphi=\{f_n\circ \varphi:n\in
\mathbb{N}\}$.
\end{lem}

\begin{pf}
It follows directly from the identity \eqref{eq17} and the fact of
$\mathcal {P}^*(\varphi^{-1}\mathcal {W})=\varphi^{-1}\mathcal
{P}^*(\mathcal {W})$ for each $\mathcal{W}\in \mathcal{C}_X$.
\end{pf}

\begin{lem}[\cite {Walters}, Lemma 9.9]\label{lem3}
Let $a_1,\cdots,a_k$ be given real numbers. If $p_i\geq0,
i=1,\cdots,k$, and $\sum_{i=1}^kp_i=1$, then
$$
\sum_{i=1}^kp_i(a_i-\log p_i)\leq \log(\sum_{i=1}^k e^{a_i}),
$$
and equality holds iff $p_i=\frac{e^{a_i}}{\sum_{i=1}^k e^{a_i}}$
for all $i=1,\cdots, k$.
\end{lem}

\begin{prop}\label{prop4}
Let $(X,T)$ be a TDS, $\mathcal{F}\in \mathcal{S}_X$  and $\mathcal
{U}\in \mathcal {C}_X^o$. Let $\pi:(X,T)\rightarrow (Y,S)$ be a
factor map between TDS, $\nu\in \mathcal {M}(Y,S)$. If $\mu \in
\mathcal{M}(X,T)$ and $ \pi\mu=\nu$, then
\begin{equation}
h_{\mu}(T,\mathcal{U}\mid Y)+ \mu(\mathcal{F})\leq \int_Y
P(T,\mathcal{F},\mathcal{U},y)d\nu(y).
\end{equation}
\end{prop}

\begin{pf}
Let $\mu=\int_Y\mu_yd\nu(y)$ be the disintegration of $\mu$ over
$\pi\mu=\nu$. As $\pi $ is a continuous map on a separable compact
space we can choose the measures $\mu_y$ such that
$\mu_y(\pi^{-1}(y))=1$ for each $y$ \cite{Bourbaki}. Then by Lemma
\ref{lem15}, we have
\begin{equation}\label{eq15}
\begin{split}
 h_{\mu}(T,\mathcal{U}\mid
Y)+\mu(\mathcal{F})=&\lim_{n\rightarrow
\infty}\frac{1}{n}H_{\mu}(\mathcal{U}_0^{n-1}\mid
Y)+\mu(\mathcal{F})\\
=&\lim_{n\rightarrow
\infty}\int_Y\frac{1}{n}H_{\mu_y}(\mathcal{U}_0^{n-1})d\nu(y)+\mu(\mathcal{F})\\
=&\lim_{n\rightarrow
\infty}\frac{1}{n}\big(\int_YH_{\mu_y}(\mathcal{U}_0^{n-1})d\nu(y)
+\int_Xf_n(x)d\mu(x) \big)\\
=&\lim_{n\rightarrow
\infty}\frac{1}{n}\int_Y\big(H_{\mu_y}(\mathcal{U}_0^{n-1})+\int_{\pi^{-1}(y)}f_n(x)d\mu_y
\big)d\nu(y).
\end{split}
\end{equation}

For any $n\in \mathbb{N}$, we have by \eqref{eq17} that there exists
a finite partition $\beta\in \mathcal {P}^*(\mathcal{U}_0^{n-1})$
such that $P_n(T,\mathcal{F},\mathcal{U},y)=
\sum\limits_{B\in\beta,B\cap \pi^{-1}(y)\neq \emptyset}
\sup\limits_{x\in B\cap \pi^{-1}(y)} \exp f_n(x)$. Let
$\beta^{'}_y=\{C: C=B\cap \pi^{-1}(y ) \,\,{\rm for \,\, some}\,\,
B\in \beta \}$, then $\beta'_y$ is a partition of $\pi^{-1}(y) $
with respect to $\beta$, and set $\beta'=\bigcup_{y\in Y}\beta'_y$.
It follows from Lemma \ref{lem3} that
\begin{equation}\label{eq16}
\begin{split}
\log(P_n(T,\mathcal{F},\mathcal{U},y))=&\log (\sum_{C\in \beta'}\sup_{x\in C}\exp f_n(x))\\
&\geq \sum_{C\in \beta'}\mu_y(C)(\sup_{x\in C}f_n(x)-\log \mu_y(C))\\
&=H_{\mu_y}(\beta')+\sum_{C\in \beta'}\sup_{x\in C}f_n(x)\cdot \mu_y(C)\\
&\geq H_{\mu_y}(\beta')+\int_{\pi^{-1}(y)}f_n(x)d\mu_y\\
&\geq H_{\mu_y}(\mathcal{U}_0^{n-1})+\int_{\pi^{-1}(y)}f_n(x)d\mu_y.
\end{split}
 \end{equation}

Combining \eqref{eq15} and \eqref{eq16}, by Fatou's Lemma and Lemma
\ref{lem2}, we have
\begin{equation}
\begin{split}
h_{\mu}(T,\mathcal{U}\mid Y)+\mu(\mathcal{F})&\leq
\limsup_{n\rightarrow \infty}\frac{1}{n}\int_Y\log
P_n(T,\mathcal{F},\mathcal{U},y)d\nu(y)\\
&\leq \int_Y \limsup_{n\rightarrow \infty}\frac{1}{n}\log
P_n(T,\mathcal{F},\mathcal{U},y)d\nu(y)\\
&=\int_YP(T,\mathcal{F},\mathcal{U},y)d\nu(y),
\end{split}
\end{equation}
and we complete the proof.
\end{pf}

The following corollary comes directly from Proposition \ref{prop4}
and the definition of $P(T,\mathcal{F},\mathcal{U}|Y)$.

\begin{cor}
Let $(X,T)$ be a TDS, $\mathcal{F}\in \mathcal{S}_X$ and $\mathcal
{U}\in \mathcal {C}_X^o$. Let $\pi:(X,T)\rightarrow (Y,S)$ be a
factor map between TDS. If $\mu \in \mathcal{M}(X,T)$, then
\begin{equation*}
h_{\mu}(T,\mathcal{U}\mid Y)+ \mu(\mathcal{F})\leq
P(T,\mathcal{F},\mathcal{U}|Y).
\end{equation*}
\end{cor}

 \begin{lem}[\cite{Xianfeng2009}, Lemma 4.4]\label{lem1}
Let $(X,T)$ be a zero-dimensional TDS. $\pi:(X,T)\rightarrow (Y,S)$
is a factor map between TDS, $y\in Y$ and $\mathcal{U}\in
\mathcal{C}^o_{X}$. Assume that for some $K\in \mathbb{N}$,
$\{\alpha_l\}_{l=1}^K$ is a sequence of finite clopen partitions of
$X$ which are finer than $\mathcal{U}$. Then for each $N\in
\mathbb{N}$, there exists a finite subset $B_N\subset \pi^{-1}(y)$
such that each atom of $(\alpha_l)_0^{N-1}, l=1,\cdots,K,$ contains
at most one point of $B_N$, and $\sum_{x\in B_N}\exp f_N(x)\geq
\frac{1}{K}P_N(T,\mathcal{F},\mathcal{U},y)$.
\end{lem}

\begin{lem}[\cite{Cao}, Lemma 2.3]\label{lem2.10}
For a sequence probability measures $\{\mu_n\}_{n=1}^{\infty}$ in
 $\mathcal{M}(X)$, where $\mu_n
 =\frac{1}{n}\sum_{i=0}^{n-1}\nu_n \circ T^{-i}$ and $\{\nu_n\}_{n=1}^{\infty}\subset
 \mathcal{M}(X)$, if $\{n_i\}$ is some
 subsequence of natural numbers $\mathbb{N}$ such that $\mu_{n_i}\rightarrow \mu \in
 \mathcal{M}(X,T)$, then for any $k\in
 \mathbb{N}$,
\begin{equation}\label{seq1}
\limsup_{i\rightarrow \infty}\frac{1}{n_i} \int f_{n_i} \,
d\nu_{n_i} \leq \frac{1}{k}\int f_k \,d\mu.
\end{equation}
In particular, the left part is no more than $\mathcal{F}_*(\mu)$.
\end{lem}

For a fixed $\mathcal{U}=\{U_1,\cdots,U_M\}\in \mathcal{C}_X^o$, we
let $\mathcal{U}^*=\{\{A_1,\cdots,A_M\}\in \mathcal{P}_X:A_m\subset
U_m, m\in\{1,\cdots,M\}\}$, where $A_m$ can be empty for some values
of $m\in\{1,\cdots,M\}$.

The following lemma will be used in the computation of
$H_{\mu}(\mathcal{U}\mid Y)$ and $h_{\mu}(T,\mathcal{U}\mid Y)$.

 \begin{lem}[\cite{huang2004}, Lemma 2]\label{lem14}
Let $G:\mathcal{P}_X\rightarrow \mathbb{R}$ be monotone in the sense
that $G(\alpha)\geq G(\beta)$ where $\alpha\succeq \beta$. Then
\begin{equation*}
\inf_{\alpha\in\mathcal{P}_X,\alpha\succeq
\mathcal{U}}G(\alpha)=\inf_{\alpha\in \mathcal{P}^*(\mathcal
{U})}G(\alpha).
\end{equation*}
\end{lem}

 \begin{prop}\label{prop2}
Let $(X,T)$ be an invertible zero-dimensional TDS, $\mathcal{F}\in
\mathcal{S}_X$ and $\mathcal{U}\in \mathcal{C}_X^o$. Let
$\pi:(X,T)\rightarrow (Y,S)$ be a factor map between TDS,
$\nu\in\mathcal{M}(Y,S)$, and $y$ be a generic point for $\nu$. Then
there exists $\mu\in\mathcal{M}(X,T)$ with $\pi\mu=\nu$ such that
\begin{equation}
P(T,\mathcal{F},\mathcal{U},y)\leq h_{\mu}^+(T,\mathcal{U}\mid
Y)+\mu(\mathcal{F}).
\end{equation}
\end{prop}

\begin{pf}
 Let $\mathcal{U}=\{U_1,U_2,\cdots,U_d\}$ and define
$$
\mathcal{U}^*=\{\alpha\in\mathcal {P}_X:
\alpha=\{A_1,A_2,\cdots,A_d\}, A_m\subset U_m, m=1,2,\cdots,d\}
$$

Since $X$ is zero-dimensional, the family of partitions in
$\mathcal{U}^*$, which are finer than $\mathcal{U}$ and consist of
clopen sets, is countable. We let $\{\alpha_l:l\geq1\}$ denote an
enumeration of this family.

Let $n\in \mathbb{N}$. By Lemma \ref{lem1}, there exists a finite
subset $B_n$ of $\pi^{-1}(y)$ such that
\begin{equation}\label{eq6}
\sum_{x\in B_n}\exp f_n(x)\geq
\frac{1}{n}P_n(T,\mathcal{F},\mathcal{U},y),
\end{equation}
and each atom of $(\alpha_l)^{n-1}_0$ contains at most one point of
$B_n$, for  all $l=1,2,\cdots,n$. Let
$$
\sigma_n=\sum_{x\in B_n}\lambda_n(x)\delta_x,
$$
where $\lambda_n(x)=\frac{\exp f_n(x)}{\sum_{y\in B_n}\exp f_n(y)}$
for $x\in B_n$, and let
$\mu_n=\frac{1}{n}\sum_{i=0}^{n-1}T^i\sigma_n$. Then
$\pi\sigma_n=\delta_y$ and
$\pi\mu_n=\frac{1}{n}\sum_{i=0}^{n-1}\delta_{S^iy}$. Choose a
subsequence $\{n_j\}$ so that $\mu_{n_j}$ converges and
$P(T,\mathcal{F},\mathcal{U},y)=\limsup_{j\rightarrow
\infty}\frac{1}{n_j}\log P_{n_j}(T,f,\mathcal{U},y)$. Let
$\mu_{n_j}\rightarrow \mu$. Then $\pi\mu=\nu$, $\mu \in
\mathcal{M}(X,T)$ and $\limsup_{n\rightarrow \infty}
\frac{1}{n_i}\int f_{n_i}d\sigma_{n_i}\leq \mu(\mathcal{F})$.

By Lemma \ref{lem14} and the fact that
\begin{equation*}
h^+_{\mu}(T,\mathcal{U}\mid Y)=\inf_{\beta\in
\mathcal{U}^*}h_{\mu}(T,\beta\mid
Y)=\inf_{l\in\mathbb{N}}h_{\mu}(T,\alpha_l\mid Y),
\end{equation*}
it is sufficient to show that for each $\l\in \mathbb{N}$,
\begin{equation*}
P(T,\mathcal{F},\mathcal{U},y)\leq h_{\mu}(T,\alpha_l\mid
Y)+\mu(\mathcal{F}).
\end{equation*}

Since $\sigma_n$ is supported on $\pi^{-1}(y)$, $T^i\sigma_n$ is
supported on $\pi^{-1}(S^iy)$ for each $i \in \mathbb{N}$, and then
$ H_{T^i\sigma_n}((\alpha_l)_0^{n-1}\mid
Y)=H_{T^i\sigma_n}((\alpha_l)_0^{n-1}) $ for each $0\leq i <n$ and
$1\leq l\leq n$.

Fix $\l\in \mathbb{N}$. For each  $n\geq l$, we know that from the
construction of $B_n$ that each atom of $(\alpha_l)_0^{n-1}$
contains at most one point in $B_n$, and,
\begin{equation}\label{eq7}
\sum_{x\in B_n}-\lambda_n(x)\log
\lambda_n(x)=H_{\sigma_n}((\alpha_l)_0^{n-1}).
\end{equation}
Combining \eqref{eq6} and \eqref{eq7}, we get that
\begin{align*}
\log P_n(T,\mathcal{F},\mathcal{U},y)-\log n &\leq \log( \sum_{x\in B_n}\exp f_n(x))\\
&=\sum_{x\in B_n}\lambda_n(x)(f_n(x)-\log \lambda_n(x))\\
&=H_{\sigma_n}((\alpha_l)_0^{n-1})+\sum_{x\in B_n}\lambda_n(x)f_n(x)\\
&=H_{\sigma_n}((\alpha_l)_0^{n-1})+\int_Xf_n(x)d\sigma_n(x).
\end{align*}
Hence
\begin{equation}\label{eq8}
\log P_n(T,\mathcal{F},\mathcal{U},y)-\log n\leq
H_{\sigma_n}((\alpha_l)_0^{n-1}\mid Y)+\int_Xf_n(x)d\sigma_n(x).
\end{equation}

Fix natural numbers $m,n$ with $n>l$ and $1\leq m\leq n-1$. Let
$a(j)=[\frac{n-j}{m}], j=0,1,\cdots,m-1$, where $[a]$ denotes the
integral part of a real number $a$. Then
\begin{equation}\label{eq9}
\bigvee_{i=0}^{n-1}T^{-i}\alpha_l=\bigvee_{r=0}^{a(j)-1}T^{-(mr+j)}(\alpha_l)^{m-1}_0\vee\bigvee_{t\in
S_j}T^{-t}\alpha_l,
\end{equation}
where $S_j=\{0,1,\cdots,j-1\}\cup\{j+ma(j),\cdots,n-1\}$. Since
${\rm card}S_j\leq 2m$, it follows from \eqref{eq8} and \eqref{eq9}
that
\begin{equation}\label{eq10}
\begin{split}
&\log P_n(T,\mathcal{F},\mathcal{U},y)-\log n\\
\leq
&\sum_{r=0}^{a(j)-1}H_{\sigma_n}(T^{-(mr+j)}(\alpha_l)^{m-1}_0\mid
Y) + H_{\sigma_n}(\bigvee_{t\in S_j}T^{-t}\alpha_l)
+\int_Xf_n(x)d\sigma_n(x)\\
\leq
&\sum_{r=0}^{a(j)-1}H_{T^{(mr+j)}\sigma_n}((\alpha_l)^{m-1}_0\mid
Y)+\int_Xf_n(x)d\sigma_n(x)+2m\log d.
\end{split}
\end{equation}
Summing up  \eqref{eq10} over $j$ from $0$ to $m-1$ then dividing
the sum by $m$ yields that
\begin{equation}\label{eq11}
\begin{split}
&\log P_n(T,\mathcal{F},\mathcal{U},y)-\log n\\
\leq
&\frac{1}{m}\sum_{j=0}^{m-1}\sum_{r=0}^{a(j)-1}H_{T^{(mr+j)}\sigma_n}((\alpha_l)^{m-1}_0\mid
Y)+\int_Xf_n(x)d\sigma_n(x)+2m\log d\\
\leq&\frac{1}{m}\sum_{j=0}^{n-1}H_{T^j\sigma_n}((\alpha_l)^{m-1}_0\mid
Y)+\int_Xf_n(x)d\sigma_n(x)+2m\log d.
\end{split}
\end{equation}
Since $H_{\{\cdot\}}((\alpha_l)^{m-1}_0\mid Y)$ is concave on
$\mathcal{M}(X)$ (Lemma 3.1 %\ref{lem12}
part (1)),
\begin{equation}\label{eq12}
\frac{1}{n}\sum_{j=0}^{n-1}H_{T^j\sigma_n}((\alpha_l)^{m-1}_0\mid
Y)\leq H_{\mu_n}((\alpha_l)^{m-1}_0\mid Y).
\end{equation}
Now by dividing \eqref{eq11} by $n$ then combining it with
\eqref{eq12}, we obtain
\begin{equation}\label{eq13}
\frac{1}{n}\log P_n(T,f,\mathcal{U},y)\leq
\frac{1}{m}H_{\mu_n}((\alpha_l)^{m-1}_0\mid
Y)+\frac{1}{n}\int_Xf_n(x)d\sigma_n(x)+\frac{2m\log d+\log n}{n}.
\end{equation}
Since $\alpha_l$ is clopen, it follows from Lemma  \ref{lem11} that
$$
\limsup_{j\rightarrow\infty}H_{\mu_{n_j}}((\alpha_l)^{m-1}_0\mid
Y)\leq H_{\mu}((\alpha_l)^{m-1}_0\mid Y).
$$
By substituting $n$ with $n_j$ in \eqref{eq13} and passing the limit
$j\rightarrow \infty$, we have that
\begin{equation}\label{eq14}
\begin{split}
P(T,\mathcal{F},\mathcal{U},y)&=\lim_{n_j\rightarrow\infty}\frac{1}{n_j}\log
P_{n_j}(T,f,\mathcal{U},y)\\
&\leq
\lim_{n_j\rightarrow\infty}\big(\frac{1}{m}H_{\mu_{n_j}}((\alpha_l)^{m-1}_0\mid
Y)+\frac{1}{n_j}\int_Xf_{n_j}(x)d\sigma_{n_j}(x)+\frac{2m\log d+\log n_j}{n_j}\big)\\
&\leq\frac{1}{m}H_{\mu}((\alpha_l)^{m-1}_0\mid Y)+\mu(\mathcal{F}).
\end{split}
\end{equation}
Then we complete the proof by taking the limit $m\rightarrow \infty$
in \eqref{eq14}.
\end{pf}

 \begin{prop}\label{prop3}
Let $(X,T)$ be an invertible zero-dimensional TDS, $\mathcal{F}\in
\mathcal{S}_X$ and $\mathcal{U}\in \mathcal{C}_X^o$. Let
$\pi:(X,T)\rightarrow (Y,S)$ be a factor map between TDS,
$\nu\in\mathcal{M}(Y,S)$. Then
\begin{equation*}
\int_YP(T,\mathcal{F},\mathcal{U},y)d\nu(y)\leq\sup\{
h^+_{\mu}(T,\mathcal{U}\mid Y)+\mu(\mathcal{F}): \mu\in
\mathcal{M}(X,T) \,\, {\rm and}\,\, \pi\mu=\nu \}.
\end{equation*}
\end{prop}

\begin{pf}
 Suppose that $\nu$ is ergodic, that is $\nu \in
\mathcal{M}^e(Y,S)$. Let $y$ be a generic point for $\nu$. By
Proposition \ref{prop2},
\begin{equation*}
P(T,\mathcal{F},\mathcal{U},y)\leq \sup_{\pi\mu=\nu}\big(
h^+_{\mu}(T,\mathcal{U}\mid Y)+\mu(\mathcal{F}) \big)=a.
\end{equation*}
Since $\nu$-{\it a.e.} $y$ is generic; so
\begin{equation*}
\int_YP(T,\mathcal{F},\mathcal{U},y)d\nu(y)\leq
\sup_{\pi\mu=\nu}\big( h^+_{\mu}(T,\mathcal{U}\mid
Y)+\mu(\mathcal{F} \big).
\end{equation*}

If $\nu $ is not ergodic, let
$\nu=\int_{\mathcal{M}^e(Y,S)}\nu_{\alpha}d\rho(\alpha)$ be its
ergodic decomposition. Let $b>0$, and
\begin{align*}
K_b=\{(\tau,\mu)\in &\mathcal{M}^e(Y,S)\times \mathcal{M}(X,T):
\pi\mu=\tau,\\
& h_{\mu}^+(T,\mathcal{U}\mid Y)+\mu(\mathcal{F}\geq \int_Y
P(T,f,\mathcal{U},y)d\tau(y)-b\}.
\end{align*}
Let $F(\tau,\mu)=F_1(\mu)-F_2(\tau)$, where
$F_1(\mu)=h_{\mu}^+(T,\mathcal{U}\mid Y)+\int_Xf(x)d\mu(x)$ and
$F_2(\tau)=\int_Y P(T,f,\mathcal{U},y)d\tau(y)$. By Lemma
\ref{lem11} and Lemma \ref{lem2}, $F_1(\mu)$ is u.s.c. on
$\mathcal{M}(X,T)$ and $F_2(\tau)$ is measurable on
$\mathcal{M}^e(Y,S)$. Moreover, $G(\mu)=F(\pi\mu,\mu)$ is measurable
on $\mathcal{M}(X,T)$. Then by the upper semi-continuity  of
$F(\tau,\cdot)$,  $F(\tau,\mu)$ is product measurable on
$\mathcal{M}^e(Y,S)\times \mathcal{M}(X,T)$. Now $K_b$ is a
measurable subset of $\mathcal{M}^e(Y,S)\times \mathcal{M}(X,T)$ and
we have shown above that $K_b$ projects onto $\mathcal{M}^e(Y,S)$.
Hence, by the selection theorem \cite{Casting}, there is a
measurable map $\phi_b:\mathcal{M}^e(Y,S)\rightarrow
\mathcal{M}(X,T)$ such that
$$
\rho(\{\tau:(\tau,\phi_b(\tau))\in K_b \})=1.
$$
Define $\mu_b$ by $\mu_b=\int_{\mathcal{M}^e(Y,S)}
\phi_b(\nu_{\alpha})d\rho(\alpha)$. Then $\mu_b\in \mathcal{M}(X,T),
\pi\mu_b=\nu$. Since $\bullet(\mathcal{F})$ is {\it u.s.c.} and
bounded affine on $\mathcal{M}(X,T)$, then by Lemma \ref{lem8} and
the well-known Choquet's Theorem (See \cite{choquet} for details),
we have
\begin{align*}
&h^+_{\mu_b}(T,\mathcal{U}\mid Y) + \mu_b(\mathcal{F})\\
=&\int_{\mathcal{M}^e(Y,S)}
h_{\phi_b(\nu(\alpha))}(T,\mathcal{U}\mid Y)d\rho(\alpha)+
\int_{\mathcal{M}^e(Y,S)}\phi_b(\nu(\alpha))(\mathcal{F})d\rho(\alpha)\\
\geq &\int_{\mathcal{M}^e(Y,S)}\big(\int_Y
P(T,\mathcal{F},\mathcal{U},y)d\nu(y)-b
\big)d\rho(\alpha)\\
=&\int_Y P(T,\mathcal{F},\mathcal{U},y)d\nu(y) -b.
\end{align*}
Therefore,
\begin{equation*}
\sup_{\pi\mu=\nu}\{ h^+_{\mu}(T,\mathcal{U}\mid
Y)+\mu(\mathcal{F})\}\geq
\int_YP(T,\mathcal{F},\mathcal{U},y)d\nu(y).
\end{equation*}
\end{pf}

 \begin{prop}\label{prop1}
Let $(X,T)$ be an invertible TDS, $\mathcal{F}\in \mathcal{S}_X$ and
$\mathcal{U}\in \mathcal{C}_X^o$. Let $\pi:(X,T)\rightarrow (Y,S)$
be a factor map between TDS, $\nu\in\mathcal{M}(Y,S)$. Then there
exists a $\mu\in \mathcal{M}(X,T)$ with $\pi\mu=\nu$ such that
\begin{equation}
h_{\mu}(T,\mathcal{U}\mid Y)+\mu(\mathcal{F})=\int_Y
P(T,\mathcal{F},\mathcal{U},y)d\nu(y).
\end{equation}
\end{prop}

\begin{pf}
 We follow the arguments in the proof of Theorem 2.5 in
\cite{Huang2006}. Let $\mathcal{U}=\{U_1,U_2,\cdots,U_M\}\in
\mathcal{C}_X^o$.

We first consider the case that $X$ is zero-dimensional, i.e., there
exists a fundamental base of the topology made of clopen sets. Since
the set of clopen subsets of $X$ is countable, the family of
partition in $\mathcal{U}^*$ consisting of clopen sets is countable.
Let $\{\alpha_l:l=1,2,\cdots\}$ ba an enumeration of this family.
Then, for any $k\in \mathbb{N}$ and $\mu \in \mathcal{M}(X,T)$, we
have
\begin{equation}\label{eq5}
h_{\mu}^+(T^k,\bigvee_{i=0}^{k-1}T^{-i}\mathcal{U}\mid
Y)=\inf_{s_k\in
\mathbb{N}^k}h_{\mu}(T^k,\bigvee_{i=0}^{k-1}T^{-i}\alpha_{s_k(i)}\mid
Y).
\end{equation}
For any $k\in \mathbb{N}$, and $s_k\in \mathbb{N}^k$, let
\begin{align*}
M(k,s_k)=\{\mu&\in
\mathcal{M}(X,T):\frac{1}{k}\big(h_{\mu}(T^k,\bigvee_{i=0}^{k-1}T^{-i}\alpha_{s_k(i)}\mid
Y) +\mu(\mathcal{F}_k)\big)\\
&\geq
\frac{1}{k}\int_YP(T^k,\mathcal{F}_k,\mathcal{U}_0^{k-1},y)d\nu(y),
\pi\mu=\nu \}.
\end{align*}
We note from Lemma \ref{lem13} that
$\frac{1}{k}\int_YP(T^k,\mathcal{F}_k,\mathcal{U}_0^{k-1},y)d\nu(y)=\int_YP(T,\mathcal{F},\mathcal{U},y)d\nu(y)$.

Since for each $k\in \mathbb{N}$, $\nu\in \mathcal{M}(Y,S^k)$, then
by Proposition \ref{prop3} there exists a $\mu_k\in
\mathcal{M}(X,T^k)$ with $\pi\mu_k=\nu$ such that
\begin{equation*}
h_{\mu_k}(T^k,\mathcal{U}_0^{k-1}\mid Y) +\mu_k(\mathcal{F}_k)\geq
\int_YP(T^k,\mathcal{F}_k,\mathcal{U}_0^{k-1},y)d\nu(y).
\end{equation*}
Since $\bigvee_{i=0}^{k-1}T^{-i}\alpha_{s_k(i)}$ is finer than
$\mathcal{U}_0^{k-1}$ for each $s_k\in \mathbb{N}^k$, we have
\begin{equation}\label{eq4}
h_{\mu}(T^k,\bigvee_{i=0}^{k-1}T^{-i}\alpha_{s_k(i)}\mid Y)
+\mu_k(\mathcal{F}_k)\geq
\int_YP(T^k,\mathcal{F}_k,\mathcal{U}_0^{k-1},y)d\nu(y).
\end{equation}

Let $\tau_k=\frac{1}{k}\sum_{i=0}^{k-1}T^i\mu_k$. Since $T^i\mu_k\in
\mathcal{M}(X,T^k), i=0,1,\cdots,k-1$, we have $\tau_k\in
\mathcal{M}(X,T)$. Moreover, since $\nu\in \mathcal{M}(Y,S)$, it is
not hard to see that $\pi\tau_k=\nu$. For $s_k\in \mathbb{N}^k$ and
$j=1,2,\cdots, k-1$, let
\begin{align*}
&P^0s_k=s_k \\
&P^js_k=\underbrace{s_k(k-j)s_k(k-j-1)\cdots s_k(k-1)}_j
\underbrace{s_k(0)s_k(1)\cdots s_k(k-1-j)}_{k-j}\in \mathbb{N}^k.
\end{align*}
It is easy to see that
\begin{align*}
&h_{T^j\mu_k}(T^k,\bigvee_{i=0}^{k-1}T^{-i}\alpha_{s_k(i)}\mid
Y)=h_{\mu_k}(T^k,\bigvee_{i=0}^{k-1}T^{-i}\alpha_{P^js_{k(i)}}\mid Y);\\
&T^j\mu_k(\mathcal{F}_k)\geq\mu_k(\mathcal{F}_k).
\end{align*}
for all $j=0,1,\cdots,k-1$. It follows from \eqref{eq4} that
\begin{equation*}
h_{T^j\mu_k}(T^k,\bigvee_{i=0}^{k-1}T^{-i}\alpha_{s_k(i)}\mid
Y)+T^j\mu_k(\mathcal{F}_k)\geq
\int_YP(T^k,\mathcal{F}_k,\mathcal{U}_0^{k-1},y)d\nu(y).
\end{equation*}

Moreover, by Lemma \ref{lem12} part(2),    for each $s_k\in
\mathbb{N}^k$,
\begin{align*}
&h_{\tau_k}(T^k,\bigvee_{i=0}^{k-1}T^{-i}\alpha_{s_k(i)}\mid
Y)+\tau_k(\mathcal{F}_k)\\
=&\frac{1}{k}\sum_{j=0}^{k-1}\big(h_{T^j\mu_k}(T^k,\bigvee_{i=0}^{k-1}T^{-i}\alpha_{s_k(i)}\mid
Y)+T^j\mu_k(\mathcal{F}_k)\big)\\
\geq &\int_YP(T^k,\mathcal{F}_k,\mathcal{U}_0^{k-1},y)d\nu(y).
\end{align*}
Hence $\tau_k\in \bigcap_{s_k\in \mathbb{N}^k}M(k,s_k)$. Let
$M(k)=\bigcap_{s_k\in \mathbb{N}^k}M(k,s_k) $. Then $M(k)$ is a
non-empty subset of $\mathcal{M}(X,T)$.

Since for every $s_k\in \mathbb{N}^k$,
$\bigvee_{i=0}^{k-1}T^{-i}\alpha_{s_k(i)}$ is a clopen cover, hence
the map
$$\mu\rightarrow h_{\mu}(T^k,
\bigvee_{i=0}^{k-1}T^{-i}\alpha_{s_k(i)}\mid Y)$$ is a u.s.c.
function from $\mathcal{M}(X,T^k)$ to $\mathbb{R}$ by Lemma
\ref{lem11} part(2). Since
$\mathcal{M}(X,T)\subset\mathcal{M}(X,T^k)$,
$h_{\{\cdot\}}(T^k,\bigvee_{i=0}^{k-1}T^{-i}\alpha_{s_k(i)}\mid Y)$
is also u.s.c. on $\mathcal{M}(X,T)$. Therefore, $M(k,s_k)$ is
closed in $\mathcal{M}(X,T)$ for each $s_k\in \mathbb{N}^k$. Thus
$M(k)$ is a non-empty closed set of $\mathcal{M}(X,T)$.

Now we show that if $k_1,k_2\in \mathbb{N}$, $k_1$ divides $k_2$,
then $M(k_2)\subset M(k_1)$. Indeed, let $\mu\in M(k_2)$ and
$k=\frac{k_2}{k_1}$. For any $s_{k_1}\in \mathbb{N}^{k_1}$, we take
$s_{k_2}=\underbrace{s_{k_1},\cdots,s_{k_1}}_k\in \mathbb{N}^{k_2}$.
Then
\begin{align*}
&\frac{1}{k_1}\big(h_{\mu}(T^{k_1},\bigvee_{i=0}^{k_1-1}T^{-i}\alpha_{s_{k_1}(i)}\mid
Y) +\mu(\mathcal{F}_{k_1})\big)\\
=&\frac{1}{k_1}\frac{1}{k}h_{\mu}(T^{kk_1},\bigvee_{j=0}^{k-1}T^{-jk_1}
\bigvee_{i=0}^{k_1-1}T^{-i}\alpha_{s_{k_1}(i)}\mid
Y)+\mu(\mathcal{F})\\
=&\frac{1}{k_2}\big(h_{\mu}(T^{k_2},\bigvee_{i=0}^{k_2-1}T^{-i}\alpha_{s_{k_2}(i)}\mid
Y) +\mu(\mathcal{F}_{k_2})\big)\\
\geq
&\frac{1}{k_2}\int_YP(T^{k_2},\mathcal{F}_{k_2},\mathcal{U}_0^{k_2-1},y)d\nu(y)\\
=&\int_YP(T,\mathcal{F},\mathcal{U},y)d\nu(y)\\
=&\frac{1}{k_1}\int_YP(T^{k_1},\mathcal{F}_{k_1},\mathcal{U}_0^{k_1-1},y)d\nu(y).
\end{align*}
Hence $\mu\in M(k_1,s_{k_1})$ for each $s_{k_1}\in \mathbb{N}^{k_1}$
and $\mu\in M(k_1)$. This shows that  $M(k_2)\subset M(k_1)$.

Since $\emptyset\neq M(k_1k_2)\subset M(k_1)\cap M(k_2)$ for any
$k_1, k_2\in \mathbb{N}$, we have that $\bigcap_{k\in
\mathbb{N}}M(k)\neq \emptyset.$

Let $\tau\in \bigcap_{k\in \mathbb{N}}M(k)$ and $k\in \mathbb{N}$,
By \eqref{eq5}, we have that
\begin{align*}
&\frac{1}{k}h_{\tau}^+(T^k,\mathcal{U}_0^{k-1}\mid
Y)+\tau(\mathcal{F})\\
=&\frac{1}{k}\big(h_{\tau}^+(T^k,\mathcal{U}_0^{k-1}\mid
Y)+k\tau(\mathcal{F}_k)   \big)\\
=&\inf_{s_k\in\mathbb{N}^k}\frac{1}{k}\big(
h_{\tau}(T^k,\bigvee_{i=0}^{k-1}T^{-i}\alpha_{s_k(i)}\mid
Y)+\tau(\mathcal{F}_k) \big)\\
\geq
&\frac{1}{k}\int_YP(T^k,\mathcal{F}_k,\mathcal{U}_0^{k-1},y)d\nu(y)=\int_YP(T,\mathcal{F},\mathcal{U},y)d\nu(y).
\end{align*}

It follows from Lemma \ref{lem10}
 that
\begin{align*}
&h_{\tau}(T,\mathcal{U}\mid
Y)+\tau(\mathcal{F})\\
=&\lim_{k\rightarrow
\infty}\frac{1}{k}\big(h_{\tau}^+(T^k,\mathcal{U}_0^{k-1}\mid
Y)+\tau(\mathcal{F}_k)   \big)\\
 \geq
&\int_YP(T,\mathcal{F},\mathcal{U},y)d\nu(y).
\end{align*}

Combining this inequality with Proposition \ref{lem9}, we complete
the proof when $X$ is zero-dimensional.

For the general case, it is well known that there exists an
invertible TDS $(Z,R)$, with $Z$ being zero-dimensional, and a
continuous surjective map $\varphi:Z\rightarrow X$ such that
$\varphi \circ R=T\circ \varphi$ (See e.g. \cite{Blanchard1997}).
For $\tau\in \mathcal{M}(Z,R)$, $\mathcal{F}\in \mathcal{S}_X$, set
$\tau (\mathcal{F}\circ \varphi)=\lim_{n\rightarrow \infty}
\frac{1}{n}\int f_n\circ \varphi d\tau$. By the above proof, we know
that there exists a $\tau\in\mathcal {M}(Z,R)$ with
$\pi(\varphi\tau)=\nu$ for the TDS $(Z,R)$ such that
$$
h_{\tau}(R,\varphi^{-1}(\mathcal{U})\mid Y)+\tau (\mathcal{F}\circ
\varphi)=\int_YP(R,\mathcal{F}\circ\varphi,
\varphi^{-1}\mathcal{U},y)d\nu(y).
$$

Let $\mu=\varphi\tau$. Then $\pi\mu=\nu$ and $\mu \in
\mathcal{M}(X,T)$. Since, by Lemma \ref{lem9},
$h_{\tau}(R,\varphi^{-1}(\mathcal{U})\mid
Y)=h_{\mu}(T,\mathcal{U}\mid Y)$, we have
\begin{equation}
\begin{split}
 &h_{\mu}(T,\mathcal{U}\mid
Y)+\mu(\mathcal{F})\\
=&h_{\tau}(R,\varphi^{-1}(\mathcal{U})\mid
Y)+\tau(\mathcal{F}\circ\varphi) =\int_YP(R,\mathcal{F}\circ\varphi,
\varphi^{-1}\mathcal{U},y)d\nu(y).
\end{split}
\end{equation}
By Lemma \ref{lem4}, we have
$$
\int_YP(R,\mathcal{F}\circ\varphi,
\varphi^{-1}\mathcal{U},y)d\nu(y)=\int_YP(T,\mathcal{F},\mathcal{U},y)d\nu(y).
$$
Then
$$
h_{\mu}(T,\mathcal{U}\mid
Y)+\mu(\mathcal{F})=\int_YP(T,\mathcal{F},\mathcal{U},y)d\nu(y),
$$
and we complete the proof of the general case.
\end{pf}

Before giving the relative local variational principle of pressure,
we first recall the notion of natural extension, which is necessary
in the proof of the relative local variational principle for the
topological pressure.

Let $d$ be the metric on $X$ and define
$\widetilde{X}=\{(x_1,x_2,\cdots): T(x_{i+1})=x_i, x_i\in X, i\in
\mathbb{N}\}$. It is clear that $\widetilde{X}$ is a subspace of the
product space $\Pi_{i=1}^{\infty}X$ with the metric $d_T$ defined by
$$
d_T((x_1,x_2,\cdots),(y_1,y_2,\cdots))=\sum_{i=1}^{\infty}\frac{d(x_i,y_i)}{2^i}.
$$
Let $\sigma_T: \widetilde{X}\rightarrow\widetilde{X}$ be the shift
homeomorphism, i.e., $\sigma_T(x_1,x_2,\cdots)=
(T(x_1),x_1,x_2,\cdots).$ We refer the TDS
$(\widetilde{X},\sigma_T)$ as the {\it natural extension} of
$(X,T)$. Let $\pi_1:\widetilde{X}\rightarrow X$ be the natural
projection onto the first component. Then
$\pi_1:(\widetilde{X},\sigma_T)\rightarrow (X,T)$ is a factor map.

Now we prove Theorem \ref{TH}, i.e., {\it let $(X,T)$ be a TDS,
$\mathcal{F}\in \mathcal{S}_X$ and $\mathcal{U}\in \mathcal{C}_X^o$,
$\pi:(X,T)\rightarrow (Y,S)$ be a factor map between TDS and $\nu\in
\mathcal{M}(Y,S)$, then
\begin{equation*}
\sup_{\mu\in\mathcal{M}(X,T)}\{h_{\mu}(T,\mathcal{U}|Y)
+\mu(\mathcal{F}):\pi\mu=\nu\}=\int_YP(T,\mathcal{F},\mathcal{U},y)d\nu(y).
\end{equation*}
}

 \begin{pf}[Proof of Theorem \ref{TH}]
 Let $(\widetilde{X},\sigma_T)$ be the natural extension of $(X,T)$
defined above. By Proposition \ref{prop1}, there exists a $\tau \in
\mathcal{M}(\widetilde{X},\sigma_T)$ such that
\begin{equation*}
h_{\tau}(\sigma_T,\pi_1^{-1}(\mathcal{U})\mid
Y)+\tau(\mathcal{F}\circ \pi_1)
=\int_YP(\sigma_T,\mathcal{F}\circ\pi_1,\pi_1^{-1}\mathcal{U},y)d\nu(y).
\end{equation*}
Let $\mu=\pi_1\tau$. Then $\mu\in\mathcal{M}(X,T)$. Since, by Lemma
\ref{lem9},
\begin{equation}\label{eq1}
h_{\mu}(T,\mathcal{U}\mid Y)+\mu(\mathcal{F})=\int_Y
P(\sigma_T,\mathcal{F}\circ\pi_1,\pi_1^{-1}\mathcal{U},y)d\nu(y).
\end{equation}
By Lemma \ref{lem4},
\begin{equation}\label{eq2}
P(\sigma_T,\mathcal{F}\circ
\pi,\pi_1^{-1}\mathcal{U},y)=P(T,\mathcal{F},\mathcal{U},y).
\end{equation}
Combining \eqref{eq1} and \eqref{eq2}, we have
\begin{equation*}\label{eq3}
h_{\mu}(T,\mathcal{U}\mid Y)+\mu(\mathcal{F})=\int_Y
P(T,\mathcal{F},\mathcal{U},y)d\nu(y).
\end{equation*}
The proof is now completed.
\end{pf}

If $(Y,S)$ is a trivial system and $\mathcal{F}=\{f\}$, then by
Lemma 2.7 in \cite{huang2007} and Theorem \ref{TH}, it is not hard
to see that Theorem \ref{TH} generalizes the standard variational
principle stated in \cite{Walters}.

Using the method to prove the outer variational principle for
entropy (\cite{Dow}), Yan {\it et al.} \cite{Yan} proved the local
outer variational principle for pressure in the single potential
case.  We shall give the following result for subadditive  sequence
of potentials without proof. For the details of the proof, we refer
the readers to see  Theorem 3 in \cite{Dow} or Theorem 2.1 in
\cite{Yan} .

\begin{lem}\label{lem17}
Let $(X,T)$ be a TDS, $\mathcal{F}\in \mathcal{S}_X$ and
$\pi:(X,T)\rightarrow (Y,S)$ be a factor map between TDS. For given
$\mathcal{U}\in \mathcal{C}_X^o$,
\begin{equation*}
P(T,\mathcal{F},\mathcal{U}|Y)=\max_{\nu\in
\mathcal{M}(Y,S)}\int_YP(T,\mathcal{F},\mathcal{U},y)d\nu(y).
\end{equation*}
\end{lem}

By Lemma \ref{lem17} and Theorem \ref{TH}, we immediately know that
Theorem \ref{coro1} holds, i.e., {\it let $(X,T)$ be a TDS,
$\mathcal{F}\in \mathcal{S}_X$, $\mathcal{U}\in \mathcal{C}_X^o$,
and $\pi:(X,T)\rightarrow (Y,S)$ be a factor map between TDS, then
\begin{equation*}
\sup\{h_{\mu}(T,\mathcal{U}\mid Y)+\mu(\mathcal{F}): \mu\in
\mathcal{M}(X,T)\}=P(T,\mathcal{F},\mathcal{U}|Y).
\end{equation*}
}

Note that for the trivial system $(Y,S)$, Theorem \ref{coro1} is
just the result obtained in \cite{zhang2009}.

\section{Pressures determine local  measure-theoretic conditional
entropies}\label{sec5}

In this section, we will prove the relative local pressure
determines the local conditional entropies.

By Theorem \ref{coro1}, it is not hard to verify that the following
results holds.
\begin{lem}
Let $(X,T)$ be a TDS,  $\mathcal{U}\in \mathcal{C}_X^o$ and
$\pi:(X,T)\rightarrow (Y,S)$ be a factor map between TDS. For any
$\mathcal{F},\mathcal{G}\in \mathcal{S}_X$ and $c\in \mathbb{R}$,
\begin{enumerate}[i)]
\item
 $P(T,\{0\},\mathcal{U}|Y)=h(T,\mathcal{U}|Y)$,
\item
If $\mathcal{F}\leq \mathcal{G}$, i.e. $f_n\leq g_n$ for all $n\in
\mathbb{N}$, then $P(T,\mathcal{F},\mathcal{U}|Y)\leq
P(T,\mathcal{G},\mathcal{U}|Y)$. In particular,
$P(T,\mathcal{F},\mathcal{U}|Y)\leq h(T,\mathcal{U}|Y)+
\|\mathcal{F}\|$,
\item
$P(T,\mathcal{F}+\{c\},\mathcal{U}|Y)=
P(T,\mathcal{F},\mathcal{U}|Y)+c$,
\item
$|P(T,\mathcal{F},\mathcal{U}|Y)-
P(T,\mathcal{G},\mathcal{U}|Y)|\leq \|\mathcal{F}-\mathcal{G}\|$,
\item
$P(T,\mathcal{\cdot},\mathcal{U}|Y) $ is convex,
\item
$P(T,\mathcal{F}+\mathcal{G}\circ T-
\mathcal{G},\mathcal{U}|Y)=P(T,\mathcal{F},\mathcal{U}|Y)$,
\item
$P(T,\mathcal{F}+\mathcal{G},\mathcal{U}|Y)\leq
P(T,\mathcal{F},\mathcal{U}|Y)+P(T,\mathcal{G},\mathcal{U}|Y)$,
\item
$P(T,c\mathcal{F},\mathcal{U}|Y)\leq c
P(T,\mathcal{F},\mathcal{U}|Y) $ if $c\geq 1$ and
$P(T,c\mathcal{F},\mathcal{U}|Y)\geq c
P(T,\mathcal{F},\mathcal{U}|Y) $ if $c\leq 1$,
\item
$|P(T,\mathcal{F},\mathcal{U}|Y)| \leq
P(T,|\mathcal{F}|,\mathcal{U}|Y) $, where
$|\mathcal{F}|=\{|f_n|:n\in \mathbb{N}\}$.
\end{enumerate}

\end{lem}

The following results shows that the relative local pressure for the
subadditive sequence of functions determines the members of
$\mathcal{M}(X,T)$. It is similar to that in the non-relative case,
and the proof can follows completely from that of Theorem 9.11 in
\cite{Walters}.

\begin{prop}
Let $\mathcal{U}\in \mathcal{C}_X^o$ and $\mu:
\mathcal{B}_X\rightarrow \mathcal{R}$ be a finite signed measure on
$X$. Then $\mu \in \mathcal{M}(X,T)$ iff $\mu(\mathcal{F})\leq
P(T,\mathcal{F},\mathcal{U}|Y)$ for all $\mathcal{F}\in
\mathcal{S}_X$.
\end{prop}

We now prove that the relative local  pressure
$P(T,\mathcal{\cdot},\mathcal{U}|Y)$ determines the local
conditional $\mu$-entropy $h_{\mu}(T,\mathcal{U}|Y)$, i.e., {\it let
$(X,T)$ be a TDS, $\mathcal{F}\in \mathcal{S}_X$ and
$\pi:(X,T)\rightarrow (Y,S)$ be a factor map between TDS, then for
given $\mathcal{U}\in \mathcal{C}_X^o$ and $\mu\in\mathcal{M}(X,T)$,
\begin{equation*}
h_{\mu}(T,\mathcal{U}|Y)=\inf\{P(T,\mathcal{F},\mathcal{U}|Y)-\mu(\mathcal{F}):\mathcal{F}\in
\mathcal{S}_X\}.
\end{equation*}
}

\begin{pf}[Proof of Theorem \ref{thm2}]
We follow the arguments in the proof of Theorem 3 in
\cite{huang2007} and Theorem 9.12 in \cite{Walters}. By Theorem
\ref{coro1}, we first have
\begin{equation*}
h_{\mu}(T,\mathcal{U}|Y)\leq
\inf\{P(T,\mathcal{F},\mathcal{U}|Y)-\mu(\mathcal{F}):\mathcal{F}\in
\mathcal{S}_X\}.
\end{equation*}

Let
\begin{equation*}
C=\{(\mu,t)\in \mathcal{M}(X,T)\times \mathbb{R}: 0\leq t\leq
h_{\mu}(T,\mathcal{U}|Y) \}.
\end{equation*}
By Theorem \ref{thm4}, the entropy map $h_{\cdot}(T,\mathcal{U}|Y):
\mathcal{M}(X,T)\rightarrow \mathbb{R}^+$ is affine. Then $C$ is
convex. Let $C(X,\mathbb{R})^*$ be the dual space of
$C(X,\mathbb{R})$ endowed with the weak*-topology and view $C$ as a
subset of $C(X,\mathbb{R})^*\times \mathbb{R}$. Take
$b>h_{\mu}(T,\mathcal{U}|Y)$. Since, by Theorem \ref{thm4}, the
entropy map $h_{\cdot}(T,\mathcal{U}|Y)$ is upper semi-continuous at
$\mu$, we have that $(\mu,b)\not \in cl (C)$. Let
$V=C(X,\mathbb{R})^*\times \mathbb{R}$, $K_1=cl(C)$,
$K_2=\{(\mu,b)\}$. Then $V$ is a locally convex, linear topological
space, and $K_1$, $K_2$ are disjoint, closed, and convex subsets of
$V$. It follows from \cite{Dun} (pp.417) that there exists a
continuous, real-valued, and convex subsets $F$ on $V$ such that
$F(x)<F(y)$ for all $x\in K_1$, $y\in K_2$, i.e.
$F:C(X,\mathbb{R})^*\times \mathbb{R}\rightarrow \mathbb{R}$ is a
continuous linear function such that $F(\mu_*,t)<F(\mu,b)$ for all
$(\mu_*,t)\in cl(C)$. Note that under the weak*-topology on $C\in
C(X,\mathbb{R})^*$, $F$ must have the form $F(\mu_*,t)=\int_X
f(x)d\mu_*(x)+td$ for some $f\in C(X,\mathbb{R})$ and some $d\in
\mathbb{R}$, i.e. $F(\mu_*,t)=\mu_*(\{f\})+td$. In particular,
$\mu_*(\{f\})+ d h_{\mu_*}(T,\mathcal{U}|Y)<\mu(\{f\})+db$
 for all $\mu_*\in \mathcal{M}(X,T)$. By taking $\mu_*=\mu$, we have
 that $d h_{\mu}(T,\mathcal{U}|Y)<db$. Hence $d>0$ and
\begin{equation*}
\mu_*(\{\frac{f}{d}\})+h_{\mu_*}(T,\mathcal{U}|Y)=\frac{\mu_*(\{f\})}{d}+
h_{\mu_*}(T,\mathcal{U}|Y)<b+\frac{\mu(\{f\})}{d}=b+\mu(\{\frac{f}{d}\}),
\end{equation*}
for all $ \mu_*\in \mathcal{M}(X,T)$. By Theorem \ref{coro1}, we
have
\begin{equation*}
P(T,\{\frac{f}{d}\},\mathcal{U}|Y)\leq b+\mu(\{\frac{f}{d}\}),
\end{equation*}
i.e.,
\begin{equation*}
b\geq P(T,\{\frac{f}{d}\},\mathcal{U}|Y)-\mu(\{\frac{f}{d}\})\geq
\inf\{
P(T,\{\mathcal{G}\},\mathcal{U}|Y)-\mu(\{\mathcal{G}\}):\mathcal{G}\in
\mathcal{S}_X \}.
\end{equation*}
Since the above inequality holds for  arbitrary $b$ satisfied
$b>h_{\mu}(T,\mathcal{U}|Y)$, we have $h_{\mu}(T,\mathcal{U}|Y)\geq
\inf\{
P(T,\{\mathcal{G}\},\mathcal{U}|Y)-\mu(\{\mathcal{G}\}):\mathcal{G}\in
\mathcal{S}_X \} $.

\end{pf}

We need the following well-known Rohlin lemma (See e.g.
\cite{Glasner}).

\begin{lem}\label{Rohlin}
 Let $(X,T)$ be invertible and $\mu \in \mathcal{M}^e(X,T)$. If $\mu $
is non-atomic, then for any $N\in \mathbb{N}$ and $\epsilon>0$,
there exists a Borel subset $D$ of $X$ such that $D,
TD,\cdots,T^{N-1}D$ are pairwise disjoint and $\mu
(\bigcup_{i=0}^{N-1}T^iD)>1-\epsilon.$

\end{lem}

We are ready to prove Theorem \ref{thm3}, i.e., {\it  let $(X,T),
(Y,S)$ be invertible TDSs, $\mathcal{F}\in \mathcal{S}_X$,
$\pi:(X,T)\rightarrow (Y,S)$ be a factor map between TDS, then for
given $\mathcal{U}\in\mathcal{C}_X^o$ and $\mu \in
\mathcal{M}(X,T)$,
\begin{equation*}
h^+_{\mu}(T,\mathcal{U}|Y)\leq
\inf\{P(T,\mathcal{F},\mathcal{U}|Y)-\mu(\mathcal{F}):\mathcal{F}\in
\mathcal{S}_X\}.
\end{equation*}
}

\begin{pf}[Proof of Theorem \ref{thm3}]
We follows the ideas in \cite{GW}, \cite{Huang2006} and
\cite{huang2007}. Since $\bullet(\mathcal{F})$ is upper
semi-continuous and bounded affine on $\mathcal{M}(X,T)$, then by
Lemma \ref{lem8} and the well-known Choquet's Theorem, it is enough
to assume that $\mu \in \mathcal{M}^e(X,T)$ and non-atomic. Then
$\nu=\pi\mu\in \mathcal{M}^e(Y,S)$. Since
$P(T,\mathcal{F}+\{c\},\mathcal{U}|Y)-\mu(\mathcal{F}+\{c\})=P(T,\mathcal{F},\mathcal{U}|Y)-\mu(\mathcal{F})$
for each $c\in \mathbb{R}$ and $\mathcal{F}\in \mathcal{S}_X$, then
we can assume that $\mathcal{F}\geq 0$, i.e. $f_n(x)\geq 0$ for each
$n\in \mathbb{N}$ and $x\in X$. Let
$\mathcal{U}=\{U_1,\cdots,U_k\}$.

For $\epsilon >0$ and $N\in \mathbb{N}$ large enough such that
\begin{equation}\label{eq18}
P_N(T,\mathcal{F},\mathcal{U},Y)\leq
2^{N(P(T,\mathcal{F},\mathcal{U}|Y)+\epsilon)} \,\, \text{and}\,\,
-(1-\frac{1}{N})\log (1-\frac{1}{N})-\frac{1}{N}\log \frac{1}{N}\leq
\epsilon.
\end{equation}
Choose small enough $1>\delta>0$ such that
\begin{equation}\label{eq25}
\sqrt{\delta}(\log k+\|f_1\| +\log (K e^{\|f_1\|}))<\epsilon.
\end{equation}
By Lemma \ref{Rohlin}, we can find a Borel subset $D$ of $X$ such
that $D,TD,\cdots, T^{N-1}D$ are pairwise disjoint and
$\mu(\bigcup_{i=0}^{N-1}T^iD)>1-\delta$. By Lemma \ref{lem16}, we
may take $\beta \in \mathcal{P}_X$ with $\beta\succeq
\mathcal{U}^{N-1}_0$ such that for each $y\in Y$,
\begin{equation}
1\leq\sum_{B\in \beta\cap \pi^{-1}(y)}\sup_{x\in B}(\exp f_N(x))\leq
P_N(T,\mathcal{F},\mathcal{U},Y).
\end{equation}
Let $\beta_D=\{B\cap D: b\in\beta\}$ be the partition of $D$. For
each $P\in \beta_D$ we can find a $s_P\in \{1,\cdots,k\}^N$ such
that $P\subset (\bigcap_{i=0}^{N-1}T^{-j}U_{i_j})\cap D$. We use the
partition $\beta_D$ to define a partition $\alpha$ of $X$ as
follows. First, for each $i=1,\cdots, k$, let
\begin{equation*}
A'_i=\bigcup_{j=0}^{N-1}\bigcup\{T^jP: P\in \beta_D\,\,
\text{and}\,\, s_P(j)=i\}.
\end{equation*}
Then let $B_1'=U_1,B_2'=U_2\backslash B_1',\cdots,
B_k'=U_k\backslash(\bigcup_{j=1}^{k-1}B_j')$. Finally, let
$A_i=A_i'\cup (B_i'\cap(X\backslash \bigcup_{j=0}^{N-1}T^jD))$ for
$i=1,\cdots,k$. Clearly, $\alpha=\{A_i:i=1,\cdots,k\}$ is a
partition of $X$ and $A_i\subset U_i$ for all $i=1,\cdots, k$. Hence
$\alpha\succeq \mathcal{U}$.

For $\beta'\in \mathcal{P}_X$ and $ R\subset X$, we define
$\beta'\cap R=\{A\cap R: A\in \beta'\,\, \text{and}\,\, A\cap R\neq
\emptyset\}$. From the construction of $\alpha$, it is easy to see
that $\alpha_0^{N-1}\cap D=\beta_D$, and moreover, for each $y\in
Y$,
\begin{align}\label{eq20}
&\sum_{C\in\alpha_0^{N-1}\cap D\cap \pi^{-1}(y)}\sup_{x\in C}(\exp
f_N(x))\notag\\
&=\sum_{C\in \beta_D\cap \pi^{-1}(y)}\sup_{x\in C}(\exp
f_N(x))\\
&\leq \sum_{C\in \beta\cap \pi^{-1}(y)}\sup_{x\in C}(\exp
f_N(x))\leq P_N(T,\mathcal{F},\mathcal{U},Y).\notag
\end{align}

Let $E=\bigcup_{i=0}^{N-1}T^iD$. Then $\mu (E)>1-\delta$. Fix $n\gg
N$, and let $G_n=\{x\in X:
\frac{1}{n}\sum_{i=0}^{n-1}1_E(T^ix)>1-\sqrt\delta \}$. Since
\begin{align*}
&\mu(G_n)+(1-\sqrt\delta)(1-\mu(G_n))\\
&\geq
\int_{G_n}\frac{1}{n}\sum_{i=0}^{n-1}1_E(T^ix)d\mu(x)+\int_{X\backslash
G_n}\frac{1}{n}\sum_{i=0}^{n-1}1_E(T^ix)d\mu(x)\\
&=\int_X \frac{1}{n}\sum_{i=0}^{n-1}1_E(T^ix)d\mu(x)\\
&=\mu(E)>1-\delta,
\end{align*}
we have
\begin{equation}
\mu(G_n)>1-\sqrt\delta.
\end{equation}

For each $x\in G_n$, let $S_n(x)=\{i\in \{0,1,\cdots,n-1\}: T^ix\in
D\}$ and $U_n(x)=\{i\in\{0,1,\cdots,n-1\}: T^ix\in E\}$. Note that
for any $x\in X$ and $i\in \mathbb{Z}$, if $T^x\in E$ then there
exists a $j\in \{0,1,\cdots,N-1\}$ such that $T^{i-j}x\in D$. Using
this fact, it is not hard to see that for each $x\in G_n$,
\begin{equation*}
U_n(x)\subseteq \bigcup_{j=0}^{N-1}(S_n(x)+j)\cup\{0,1,\cdots,N-1\}.
\end{equation*}
Since for each $x\in G_n$,
$|U_n(x)|=\sum_{i=0}^{n-1}1_E()T^ix>1-\sqrt\delta$, we have
$|\{0,1,\cdots,n-1\}\backslash U_N(x)|\leq n\sqrt \delta$.
Therefore, for each $x\in G_n$,
\begin{align}\label{eq19}
&|\{0,1,\cdots,n-1\}\backslash \bigcup_{j=0}^{N-1}(S_n(x)+j)| \notag\\
\leq &|\{0,1,\cdots,N-1\}\cup\{0,1,\cdots,n-1\}\backslash U_n(x)|\notag\\
\leq &n\sqrt\delta +N.
\end{align}

Let $\mathcal{F}_n=\{S_n(x):x\in G_n\}$. Since for each $F\in
\mathcal{F}_n$, $F\cap(F+i) =\emptyset, i=1,\cdots,N-1$, we have
$|F|\leq \frac{n}{N}+1$. Hence
\begin{equation*}
|\mathcal{F}_n|\leq \sum_{j=1}^{a_n}\frac{n!}{j!\cdot(n-j)!}\leq
a_n\frac{n!}{a_n!\cdot(n-a_n)!}\leq n\frac{n!}{a_n!\cdot (n-a_n)!}
\end{equation*}
where $a_n=[\frac{n}{N}]+1$. By Stirling's formulation and the
second inequality in \eqref{eq18}, we have
\begin{equation*}
\lim_{n\rightarrow\infty}\frac{1}{n}\log(n\frac{n!}{a_n!\cdot(n-a_n)!})=-(1-\frac{1}{N})\log
(1-\frac{1}{N})-\frac{1}{N}\log\frac{1}{N} <\epsilon.
\end{equation*}
Hence we have
\begin{equation}\label{eq26}
\limsup_{n\rightarrow \infty}\frac{1}{n}\log |\mathcal{F}_n|\leq
\lim_{n\rightarrow \infty}\frac{1}{n}\log n\frac{n!}{a_n!\cdot
(n-a_n)!}\leq \epsilon.
\end{equation}
For each $F\in \mathcal{F}_n$, let $B_F=\{x\in G_n:S_n(x)=F\}$.
Clearly, $\{B_F\}_{F\in \mathcal{F}_n}$ forms a partition of $G_n$.

For each $F\in \mathcal{F}_n$, $F=\{s_1<s_2<\cdots<s_l\}$, let
$H_F=\{0,1,\cdots,n-1\}\backslash \bigcup_{i=0}^{N-1}(F+i)$. It
follows from \eqref{eq19} that $l\leq \frac{n}{N}+1$, $|H_F|\leq
n\sqrt\delta +N$. Moreover, for each $y\in Y$, using \eqref{eq20}
and the facts that $|\alpha|=k$,
$P_N(T,\mathcal{F},\mathcal{U},Y)\geq 1$, $B_F\subseteq G_n\cap
\bigcap_{j=1}^lT^{-s_j}D$ and $f_n(x)\leq
\sum_{j=1}^lf_N(T^{s_j}x)+\sum_{r\in H_F}f_1(T^rx)$, we have
\begin{align*}
&\sum_{C\in\alpha_0^{n-1}\cap B_F\cap \pi^{-1}(y)}\sup_{x\in C}(\exp
f_n(x))\notag\\
&\leq \sum_{C\in\alpha_0^{n-1}\cap \bigcap_{j=1}^lT^{-s_j}D \cap
\pi^{-1}(y)}\sup_{x\in C}(\exp f_n(x))\notag\\
&=\sum_{C\in \bigvee_{j=1}^lT^{-s_j}\alpha_0^{N-1}\cap
\bigcap_{j=1}^lT^{-s_j}D\cap \pi^{-1}(y)\vee \bigvee_{r\in
H_F}T^{-r}\alpha}\sup_{x\in C}(\exp f_n(x))\notag\\
&=\sum_{C\in \bigvee_{j=1}^lT^{-s_j}(\alpha_0^{N-1}\cap D)\cap
\pi^{-1}(y)\vee \bigvee_{r\in H_F}T^{-r}\alpha}\sup_{x\in C}(\exp
f_n(x))
\notag\\
&\leq \sum_{C\in \bigvee_{j=1}^lT^{-s_j}(\alpha_0^{N-1}\cap D)\cap
\pi^{-1}(y)\vee \bigvee_{r\in H_F}T^{-r}\alpha}\sup_{x\in
C}(\exp(\sum_{j=1}^lf_N(T^{s_j}x)\notag\\
&\quad\quad \quad\quad \quad\quad \quad\quad \quad\quad \quad\quad
\quad\quad \quad\quad \quad \quad\quad \quad\quad \quad\quad \quad
\quad
+\sum_{r\in H_F}f_1(T^rx)))\notag\\
& \leq \sum_{C\in \bigvee_{j=1}^lT^{-s_j}(\alpha_0^{N-1}\cap D)\cap
\pi^{-1}(y)}\sup_{x\in C}(\exp(\sum_{j=1}^lf_N(T^{s_j}x)))\notag\\
 &\quad\quad \quad\quad \quad\quad \quad\quad \quad\quad \quad\quad \quad\quad \quad
 \cdot \sum_{C\in
\bigvee_{r\in H_F}T^{-r}\alpha}\sup_{x\in C}(\exp (\sum_{r\in
H_F}f_1(T^rx)))\notag\\
& \leq \prod_{j=1}^l\big(\sum_{C\in T^{-s_j}(\alpha^{N-1}_0\cap
D)\cap
\pi^{-1}(y)}\sup_{x\in C}(\exp(f_N(T^{s_j}x))) \big)\notag\\
 &\quad\quad \quad\quad \quad\quad \quad\quad \quad\quad \quad\quad \quad\quad \quad\quad
\cdot \prod_{r\in H_F}\big(\sum_{C\in T^{-r}\alpha}\sup_{x\in C}
(\exp f_1(T^rx))\big)\notag\\
&=\prod_{j=1}^l\big(  \sum_{C\in \alpha^{N-1}_0\cap D\cap
\pi^{-1}(S^{s_j}(y))}\sup_{x\in C}(\exp(f_N(x))) \big)\cdot
\big(\sum_{C\in \alpha}\sup_{x\in C} (\exp
f_1(x))\big)^{|H_F|}\notag\\
& \leq (P_N(T,\mathcal{F},\mathcal{U},Y))^{l}\cdot(k\cdot
e^{\|f_1\|})^{|H_F|} \text{\quad (by \eqref{eq20}} )     \notag\\
&\leq (P_N(T,\mathcal{F},\mathcal{U},Y))^{\frac{n}{N}+1}\cdot(k\cdot
e^{\|f_1\|})^{n\sqrt\delta +N}\notag\\
\end{align*}
Summing this result over all $F\in \mathcal{F}_n$ yields that
\begin{equation}\label{eq23}
\begin{split}
\sum_{F\in \mathcal{F}_n}\sum_{C\in\alpha_0^{n-1}\cap B_F\cap
\pi^{-1}(y)}&\sup_{x\in C}(\exp f_n(x))\\
&\leq |\mathcal{F}_n|\cdot
P_N(T,\mathcal{F},\mathcal{U},Y))^{\frac{n}{N}+1}\cdot(k\cdot
e^{\|f_1\|})^{n\sqrt\delta +N}.
\end{split}
\end{equation}

Let $\mu=\int_Y\mu_yd\nu(y)$ be the disintegration of $\mu$ over
$\pi\mu=\nu$. Choose the measures $\mu_y\in \mathcal{M}(X)$ such
that $\mu_y(\pi^{-1}(y))=1$ for each $y\in Y$. For each $F\in
\mathcal{F}_n$, we have
\begin{equation}\label{eq21}
\begin{split}
&H_{\mu}(\alpha_0^{n-1}\cap B_F|Y)+\int_{B_F}f_nd\mu\\
=&\int_YH_{\mu_y}(\alpha_0^{n-1}\cap
B_F)d\nu(y)+\int_Y\int_{B_F}f_nd\mu_yd\nu(y)\\
=&\int\big( H_{\mu_y}(\alpha_0^{n-1}\cap B_F)+ \int_{B_F}f_nd\mu_y
\big)\\
\leq &\int_Y \sum_{C\in \alpha_0^{n-1}\cap
B_F\cap\pi^{-1}(y)}\mu_y(C)(\sup_{x\in C}f_n(x)-\log\mu_y(C) )
d\nu(y).
\end{split}
\end{equation}

Since $\mu(X\backslash G_n)<\sqrt\delta$ and $|\alpha_0^{n-1}\cap
(X\backslash G_n)|\leq k^n$, we have
\begin{equation}\label{eq22}
\begin{split}
&H_{\mu}(\alpha_0^{n-1}\cap (X\backslash G_n)|Y)+\int_{X\backslash
G_n}f_nd\mu\\
=&\int_Y H_{\mu_y}(\alpha_0^{n-1}\cap (X\backslash
G_n))d\nu(y)+\int_Y \int_{X\backslash G_n}f_nd\mu_y d\nu(y)\\
=& \int_Y\big( H_{\mu_y}(\alpha_0^{n-1}\cap (X\backslash
G_n))+\int_{X\backslash G_n}f_nd\mu_y \big)\\
\leq &\int_Y\big(\sum_{C'\in \alpha_0^{n-1}\cap (X\backslash
G_n)}-\mu_y(C')\log \mu_y(C') +\mu_y(X\backslash G_n)\cdot
\|f_n\|\big)d\nu(y)\\
\leq &\int_y -\big(\sum_{C'\in \alpha_0^{n-1}\cap (X\backslash
G_n)}\mu_y(C') \big)\log \frac{\sum_{C'\in \alpha_0^{n-1}\cap
(X\backslash G_n)}\mu_y(C')}{|\alpha_0^{n-1}\cap (X\backslash
G_n)|}\\
&\quad \quad\quad \quad\quad \quad\quad \quad\quad \quad\quad
\quad\quad \quad\quad \quad\quad \quad\quad +\mu_y(X\backslash
G_n)\cdot \|f_n\|\big)d\nu(y)\\
=&\int_Y \big(-\mu_y(X\backslash G_n)\log \mu_y(X\backslash G_n) \\
&   \quad\quad \quad\quad \quad\quad \quad\quad + \mu_y(X\backslash
G_n) \big( \log |\alpha_0^{n-1}\cap (X\backslash
G_n)|+\|f_n\|\big)\big)d\nu(y)\\
\leq &\int_Y -\mu_y(X\backslash G_n)\log \mu_y(X\backslash
G_n)d\nu(y)+\sqrt\delta(\log k^n +\|f_n\|)
\end{split}
\end{equation}

Let $\gamma=\{B_F\}_{F\in \mathcal{F}_n}\cup \{X\backslash G_n\}$.
Then by \eqref{eq23}, \eqref{eq21}, \eqref{eq22} and Lemma
\ref{lem3}, we have
\begin{equation}\label{eq24}
\begin{split}
&H_{\mu}(\alpha_0^{n-1}|Y)+\int_Xf_nd\mu
\leq H_{\mu}(\alpha_0^{n-1}\vee \gamma|Y)+\int_Xf_nd\mu\\
=&\sum_{F\in \mathcal{F}_n}\big(H_{\mu}(\alpha_0^{n-1}\cap
B_F|Y)+\int_{B_F}f_nd\mu\big) \\
&\quad \quad\quad \quad\quad \quad\quad \quad\quad + \big(
H_{\mu}(\alpha_0^{n-1}\cap (X\backslash G_n)|Y)+\int_{X\backslash
G_n}f_nd\mu \big)\\
\leq &\int_Y \big( \sum_{F\in \mathcal{F}_n}\sum_{C\in
\alpha_0^{n-1}\cap
B_F\cap\pi^{-1}(y)}\mu_y(C)(\sup_{x\in C}f_n(x) -\log\mu_y(C) )\\
&\quad \quad\quad \quad\quad \quad
 -\mu_y(X\backslash G_n)(0-\log
\mu_y(X\backslash G_n)) \big)d\nu(y) +\sqrt\delta(\log k^n
+\|f_n\|)\\
\leq &\int_Y  \log\big( \sum_{F\in \mathcal{F}_n}\sum_{C\in
\alpha_0^{n-1}\cap B_F\cap\pi^{-1}(y)} e^{\sup_{x\in C}f_n(x)}
+e^{\sup_{x\in X\backslash G_n}0} \big) d\nu(y)\\
&\quad \quad\quad \quad\quad \quad \quad \quad\quad \quad\quad \quad
\quad\quad \quad\quad \quad \quad \quad\quad \quad\quad \quad
+n\sqrt\delta(\log k +\frac{\|f_n\|}{n})\\
\leq & n(b_n+\sqrt\delta(\log k+\|f_1\|)),
\end{split}
\end{equation}
where $b_n=\frac{1}{n}\log\big(|\mathcal{F}_n|\cdot
P_N(T,\mathcal{F},\mathcal{U},Y))^{\frac{n}{N}+1}\cdot(k\cdot
e^{\|f_1\|})^{n\sqrt\delta +N}+1\big)$.

Hence, by \eqref{eq18}, \eqref{eq25}, \eqref{eq26} and \eqref{eq24},
we have
\begin{equation*}
\begin{split}
&h^+_{\mu}(T,\mathcal{U}|Y)+\mu(\mathcal{F})\leq
h_{\mu}(T,\alpha|Y)+\mu(\mathcal{F})\\
=&\lim_{n\rightarrow\infty}\frac{1}{n}(H_{\mu}(\alpha_0^{n-1}|Y)+\int_Xf_nd\mu)
\leq \limsup_{n\rightarrow \infty}b_n + \sqrt\delta(\log k+
\|f_1\|)\\
=&\limsup_{n\rightarrow
\infty}\frac{1}{n}\big(\log|\mathcal{F}_n|+({\frac{n}{N}+1})\log
P_N(T,\mathcal{F},\mathcal{U},Y))\\
&\quad \quad\quad \quad\quad
 \quad\quad \quad\quad \quad \quad
 +({n\sqrt\delta +N})\log(k\cdot
e^{\|f_1\|})\big) + \sqrt\delta(\log k+
\|f_1\|)\\
=&\limsup_{n\rightarrow
\infty}\frac{1}{n}\log|\mathcal{F}_n|+\frac{1}{N}P_N(T,\mathcal{F},\mathcal{U},Y)
+\sqrt\delta(\log k+ \|f_1\| +\log(k\cdot e^{\|f_1\|}))\\
\leq &\frac{1}{N}P_N(T,\mathcal{F},\mathcal{U},Y)+2\epsilon\\
\leq &P(T,\mathcal{F},\mathcal{U}|Y)+3\epsilon.
\end{split}
\end{equation*}
Since $\epsilon>0$ is arbitrary, then the proof of Theorem
\ref{thm3} is completed.
\end{pf}

For $\mathcal{F}=\{0\}$, by Theorem \ref{thm2} and \ref{thm3}, we
have $h_{\mu}^+(T,\mathcal{U}|Y)=h_{\mu}(T,\mathcal{U}|Y)$ for the
invertible TDS. Moreover, if $(Y,S)$ is the trivial system, then
$h_{\mu}^+(T,\mathcal{U})=h_{\mu}(T,\mathcal{U})$. These results
were  shown in \cite{GW}, \cite{Huang2006} and \cite{huang2007}.

\section*{Acknowledgements}
The first and second authors are  supported by the National Natural
Science Foundation of China (Grant No. 10971100). The second author
is  partially supported by National Basic Research Program of China
(973 Program) (Grant No. 2007CB814800)

\end{document}